\documentclass[11pt]{amsart}

\usepackage{amsfonts}
\usepackage{mathrsfs}
\usepackage{bbm}
\usepackage{amsmath}
\usepackage{amssymb}
\usepackage{amscd}
\usepackage{amsthm}

\newtheorem{thm}{Theorem}[section]
\newtheorem{lemma}[thm]{Lemma}
\newtheorem{prop}[thm]{Proposition}
\newtheorem{cor}[thm]{Corollary}

\newtheorem{question}{\textbf{Question}}

\newcommand{\di}{\displaystyle}

\newcommand{\ve}{\varepsilon}

\newcommand{\bc}{\mathbb{C}}
\newcommand{\br}{\mathbb{R}}
\newcommand{\bz}{\mathbb{Z}}

\newcommand{\cu}{\mathscr{C}}
\newcommand{\fc}{\mathscr{C}}
\newcommand{\sg}{\mathscr{G}}
\newcommand{\ci}{\mathscr{I}}

\newcommand{\sm}{\mathscr{M}}
\newcommand{\sq}{\mathscr{Q}}
\newcommand{\cs}{\mathscr{S}}
\newcommand{\st}{\mathscr{T}}
\newcommand{\sv}{\mathscr{V}}
\newcommand{\su}{\mathscr{U}}

\newcommand{\cb}{\mathcal{B}}

\newcommand{\cf}{\mathcal{F}}
\newcommand{\cg}{\mathcal{G}}
\newcommand{\cl}{\mathcal{L}}
\newcommand{\co}{\mathcal{O}}

\newcommand{\fg}{\mathfrak{g}}

\newcommand{\fin}{\mathfrak{P}}

\newcommand{\fO}{\mathfrak{O}}
\newcommand{\fv}{\mathfrak{V}}

\newcommand{\fo}{\mathcal{F}_1}
\newcommand{\ft}{\mathcal{F}_2}

\newcommand{\dist}{\mathrm{dist}}

\newcommand{\bp}{\mathbb{P}}
\newcommand{\py}{\mathbb{P}^1}
\newcommand{\pn}{\mathbb{P}^n}

\newcommand{\tv}{V^{\ast}}
\newcommand{\ctv}{\overline{V^{\ast}}}
\newcommand{\tw}{W^{\ast}}
\newcommand{\ctw}{\overline{W^{\ast}}}

\newcommand{\pa}{\phi^{\ast}}

\newcommand{\paotm}{\phi^{\ast}_{M}}
\newcommand{\paotn}{\phi^{\ast}_{N}}
\newcommand{\paotq}{\phi^{\ast}_{Q}}
\newcommand{\paotpn}{\phi^{\ast}_{\mathbb{P}^n}}

\newcommand{\paotrm}{\phi^{\ast}_{\mathbb{R}^m}}

\newcommand{\tovtw}{\mathscr{T}\left(\mathcal{F}_1, V, \mathcal{F}_2, W \right)}
\newcommand{\tovtwz}{\mathscr{T}\left(\mathcal{F}_1, V, \mathcal{F}_2, \overline{W_0} \right)}

\newcommand{\tfvfv}{\mathscr{T}\left(\mathcal{F}, V, \mathcal{F}, V \right)}
\newcommand{\tfv}{\mathscr{T}\left(\mathcal{F}, V\right)}

\newcommand{\tkovktw}{\mathscr{T}\left(C^{k_1, \alpha_1}, V, C^{k_2, \alpha_2}, W \right)}

\newcommand{\dfv}{\mathscr{D}\left(\mathcal{F}, V \right)}
\newcommand{\dciv}{\mathscr{D}\left(C^{\infty}, V \right)}
\newcommand{\dciso}{\mathscr{D}\left(C^{\infty}, S^1 \right)}

\newcommand{\tcisotw}{\mathscr{T}\left(C^{\infty}, S^1, \mathcal{F}_2, W \right)}

\newcommand{\auto}{\mathrm{Aut}\left(\fvmz \right)}

\newcommand{\autopn}{\mathrm{Aut}\left(\fvpnz \right)}

\newcommand{\supp}{\mathrm{supp}}

\newcommand{\df}{\mathrm{deg} f}

\newcommand{\fvm}{\mathcal{F}(V, M)}

\newcommand{\fovm}{\mathcal{F}_1\left(V, M \right)}
\newcommand{\footvm}{\mathcal{F}_1\left(\overline{V^{\ast}}, M \right)}

\newcommand{\ftwm}{\mathcal{F}_2\left(W, M \right)}

\newcommand{\fsom}{\mathcal{F}(S^1, M)}

\newcommand{\fvq}{\mathcal{F}(V, Q)}
\newcommand{\fovq}{\mathcal{F}_1\left(V, Q \right)}
\newcommand{\ftwq}{\mathcal{F}_2\left(W, Q \right)}

\newcommand{\fvn}{\mathcal{F}(V, N)}

\newcommand{\fovn}{\mathcal{F}_1\left(V, N \right)}
\newcommand{\fovnz}{\mathcal{F}^0_1\left(V, N \right)}
\newcommand{\ftwn}{\mathcal{F}_2\left(W, N \right)}

\newcommand{\fvmz}{\mathcal{F}^{0}\left(V, M \right)}

\newcommand{\fvg}{\mathcal{F}\left(V, \cg \right)}

\newcommand{\fvgz}{\mathcal{F}^0\left(V, \cg \right)}

\newcommand{\fovg}{\mathcal{F}_1\left(V, \cg \right)}
\newcommand{\ftwg}{\mathcal{F}_2\left(W, \cg \right)}
\newcommand{\cwg}{C\left(W, \cg \right)}

\newcommand{\fovgo}{\mathcal{F}^0_1\left(V, \cg \setminus \{1\} \right)}

\newcommand{\ftwgo}{\mathcal{F}_2\left(W, \cg \setminus \{1\} \right)}

\newcommand{\fovcno}{\mathcal{F}_1\left(V, \bc^n \setminus \{\zero\} \right)}

\newcommand{\ftwag}{\mathcal{F}_2\left(W, Aut\left(\cg\right) \right)}

\newcommand{\fvfg}{\mathcal{F}\left(V, \mathfrak{g} \right)}

\newcommand{\fovfg}{\mathcal{F}_1\left(V, \mathfrak{g} \right)}

\newcommand{\fovgz}{\mathcal{F}_1^0\left(V, \cg \right)}
\newcommand{\ftwgz}{\mathcal{F}_2^0\left(W, \cg \right)}

\newcommand{\fovcn}{\mathcal{F}_1\left(V, \mathbb{C}^n \right)}
\newcommand{\ftwcn}{\mathcal{F}_2\left(W, \mathbb{C}^n \right)}

\newcommand{\ftwcm}{\mathcal{F}_2\left(W, \mathbb{C}^m \right)}

\newcommand{\fovck}{\mathcal{F}_1\left(V, \mathbb{C}^k \right)}

\newcommand{\fovc}{\mathcal{F}_1\left(V, \mathbb{C} \right)}

\newcommand{\civm}{C^{\infty}\left(V, M\right)}

\newcommand{\fsopn}{\cf\left(S^1, \mathbb{P}^n \right)}
\newcommand{\fsopo}{\cf\left(S^1, \mathbb{P}^1 \right)}

\newcommand{\cisopo}{C^{\infty}\left(S^1, \mathbb{P}^1 \right)}
\newcommand{\cisopoa}{C^{\infty}\left(S^1, \mathbb{P}^1 \setminus \{a\} \right)}

\newcommand{\fosopo}{\mathcal{F}_1\left(S^1, \mathbb{P}^1 \right)}

\newcommand{\cisopm}{C^{\infty}\left(S^1, \mathbb{P}^m \right)}
\newcommand{\cisopmj}{C^{\infty}\left(S^1, \mathbb{P}^{m_j} \right)}

\newcommand{\cisoc}{C^{\infty}\left(S^1, \mathbb{C} \right)}

\newcommand{\fvpn}{\mathcal{F}\left(V, \mathbb{P}^n \right)}
\newcommand{\fvpnz}{\mathcal{F}^{0} \left(V, \mathbb{P}^n \right)}

\newcommand{\fvpy}{\mathcal{F}\left(V, \mathbb{P}^1 \right)}
\newcommand{\fvpyz}{\mathcal{F}^{0} \left(V, \mathbb{P}^1 \right)}

\newcommand{\fvpm}{\mathcal{F}\left(V, \mathbb{P}^m \right)}

\newcommand{\fovpnz}{\mathcal{F}_1^{0} \left(V, \mathbb{P}^n \right)}

\newcommand{\ftwpn}{\mathcal{F}_2 \left(W, \mathbb{P}^n \right)}

\newcommand{\fovpmj}{\mathcal{F}_1 \left(V, \mathbb{P}^{m_j} \right)}

\newcommand{\ftwpm}{\mathcal{F}_2 \left(W, \mathbb{P}^m \right)}

\newcommand{\ftwpmo}{\mathcal{F}_2 \left(W, \mathbb{P}^{m_1} \right)}
\newcommand{\ftwpmp}{\mathcal{F}_2 \left(W, \mathbb{P}^{m_p} \right)}

\newcommand{\ftwzpm}{\mathcal{F}_2 \left(\overline{W_0}, \mathbb{P}^m \right)}

\newcommand{\civpn}{C^{\infty}\left(V, \mathbb{P}^n \right)}

\newcommand{\kp}{W^{k, p}}

\newcommand{\ka}{C^{k, \alpha}}
\newcommand{\kavm}{C^{k, \alpha}\left(V, M\right)}

\newcommand{\fvcn}{\mathcal{F}\left(V, \mathbb{C}^n \right)}
\newcommand{\fvcm}{\mathcal{F}\left(V, \mathbb{C}^m \right)}

\newcommand{\cipo}{C^{\infty}\left(S^1, \mathbb{P}^1 \right)}

\newcommand{\fvrn}{\mathcal{F}\left(V, \mathbb{R}^n \right)}

\newcommand{\fvrm}{\mathcal{F}\left(V, \mathbb{R}^m \right)}

\newcommand{\cotvcm}{C^{1}(\tv, \mathbb{C}^m)}

\newcommand{\fovrn}{\mathcal{F}_1(V, \mathbb{R}^n)}

\newcommand{\fovr}{\mathcal{F}_1(V, \mathbb{R})}

\newcommand{\ftwrn}{\mathcal{F}_2(W, \mathbb{R}^n)}

\newcommand{\ftwr}{\mathcal{F}_2(W, \mathbb{R})}

\newcommand{\txly}{\tilde{y}+\lambda \tilde{x}}

\newcommand{\pxly}{\left(\tilde{y}+\lambda \tilde{x}\right)}

\newcommand{\xl}{x_{\lambda}}

\newcommand{\xlz}{x_{\lambda_0}}

\newcommand{\fytm}{\mathcal{F}\left(y^{\ast}TM\right)}
\newcommand{\foytm}{\mathcal{F}_1\left(y^{\ast}TM\right)}

\newcommand{\fxtpn}{\mathcal{F}\left(x^{\ast}T\mathbb{P}^n\right)}

\newcommand{\fgxtpm}{\mathcal{F}\left(\gamma(x)^{\ast}T\mathbb{P}^m\right)}

\newcommand{\kaytm}{C^{k, \alpha}\left(y^{\ast}TM\right)}

\newcommand{\ciytm}{C^{\infty}\left(y^{\ast}TM\right)}

\newcommand{\tx}{\tilde{x}}
\newcommand{\ty}{\tilde{y}}

\newcommand{\hax}{\hat{x}}

\newcommand{\wf}{\mathscr{W}_f}

\newcommand{\htv}{H^2(V, \mathbb{Z})}

\newcommand{\dimv}{\dim_{\br} V}
\newcommand{\dimm}{\dim_{\bc} M}

\newcommand{\yxt}{y_{\xi, t}}
\newcommand{\yxoto}{y_{\xi_1, t_1}}

\newcommand{\bv}{\beta_{\ast}^{\mathcal{F}, V}}
\newcommand{\bla}{\beta_{\ast}}

\newcommand{\pov}{\Phi_{\ast}^{\mathcal{F}_1, V}}

\newcommand{\iptw}{\left(\Phi^{-1} \right)_{\ast}^{\mathcal{F}_2, W}}

\newcommand{\hd}{\mathscr{H}_d^{n, m}}

\newcommand{\ho}{\mathscr{H}_1^{n, n}}

\newcommand{\hyn}{\mathscr{H}_1^{1, n}}

\newcommand{\hdo}{\mathscr{H}_0^{n, m}}

\newcommand{\fvhd}{\mathcal{F}\left(V, \mathscr{H}_d^{n, m} \right)}

\newcommand{\ftwhd}{\mathcal{F}_2\left(W, \mathscr{H}_d^{n, m} \right)}

\newcommand{\zero}{\mathbf{0}}

\newcommand{\fh}{\hat{f}}
\newcommand{\xh}{\hat{x}}

\newcommand{\cih}{\hat{\mathscr{I}}}

\newcommand{\vb}{\bar{v}}

\newcommand{\gf}{\gamma_f}

\newcommand{\cxy}{C_{x, \, y}}

\newcommand{\xtm}{x^{\ast} TM}
\newcommand{\xtn}{x^{\ast} TN}

\newcommand{\fvpglnp}{\mathcal{F}\left(V, PGL(n+1, \bc) \right)}
\newcommand{\fvpglt}{\mathcal{F}\left(V, PGL(2, \bc) \right)}

\newcommand{\fvpglnptr}{\mathcal{F}^I\left(V, PGL(n+1, \bc) \right)}

\newcommand{\tf}{\tilde{f}}

\newcommand{\lpo}{\Lambda_{\varphi_1}}

\newcommand{\hocipo}{\di H^0\left(C^{\infty}(S^1, \mathbb{P}^1), \Lambda_{\varphi_1}\right)}

\newcommand{\hopono}{\di H^0\left(\bp^1, \co(n_1) \right)}

\newcommand{\ahopono}{\di E_{t_1}^{\ast} H^0\left(\bp^1, \co(n_1) \right)}

\newcommand{\ehopono}{\di E^{\ast} H^0\left(\bp^1, \co(n_1) \right)}

\newcommand{\jahojk}{\di (j^k_{t_1})^{\ast} H^0\left(J^k(t_1, \bp^1), \tau_k^{\ast} \co(n_1) \right)}
\newcommand{\jahojkm}{\di (j^{k-1}_{t_1})^{\ast} H^0\left(J^{k-1}(t_1, \bp^1), \tau_{k-1}^{\ast} \co(n_1) \right)}
\newcommand{\jahojno}{\di (j^{n_1-1}_{t_1})^{\ast} H^0\left(J^{n_1-1}(t_1, \bp^1), \tau_{n_1-1}^{\ast} \co(n_1) \right)}

\newcommand{\jk}{\mathfrak{J}^k}
\newcommand{\jnom}{\mathfrak{J}^{n_1-1}}

\newcommand{\ts}{\tilde{\sigma}}

\newcommand{\fepo}{\mathfrak{E}_{\varphi_1}}

\newcommand{\cistpo}{C^{\infty}(S^2, \bp^1)}

\begin{document}

\title{Range decreasing group homomorphisms and holomorphic maps between generalized loop spaces}

\author{Ning Zhang}


\thanks{The author is grateful to L. Lempert for stimulating discussions and very helpful comments on the present paper.
This research was partially supported by the Scientific Research Foundation
of Ocean University of China grant 861701013110, and the National
Natural Science Foundation of China grants 10871002 and 11371035.}

\begin{abstract}
  Let $\cg$ resp. $M$ be a positive dimensional Lie group resp. connected complex manifold without boundary
  and $V$ a finite dimensional $C^{\infty}$ compact connected manifold, possibly with boundary.
  Fix a smoothness class $\cf=C^{\infty}$, H\"older $\ka$ or Sobolev $\kp$.
  The space $\fvg$ resp. $\fvm$ of all $\cf$ maps $V \to \cg$ resp. $V \to M$
  is a Banach/Fr\'echet Lie group resp. complex manifold.
  Let $\fvgz$ resp. $\di \fvmz$ be the component of $\fvg$ resp. $\fvm$ containing the identity resp. constants.
  A map $f$ from a
  domain $\Omega \subset \fovm$ to $\ftwm$ is called range decreasing if $f(x)(W) \subset x(V)$, $x \in \Omega$.
  We prove that if $\dim_{\br} \cg \ge 2$, then any range decreasing group homomorphism $f: \fovgz \to \ftwg$
  is the pullback by a map $\phi: W \to V$.
  We also provide several sufficient conditions for a range decreasing holomorphic map
  $\Omega$ $\to$ $\ftwm$ to be a pullback operator. Then we apply these results to study certain decomposition of
  holomorphic maps $\di \fovn \supset \Omega \to \ftwm$.
  In particular, we identify some classes of holomorphic maps $\di \fovpnz \to \ftwpm$,
  including all automorphisms of $\fvpnz$.
\end{abstract}


\address{School of Mathematical Sciences \\ Ocean University of China \\Qingdao, 266100, P. R. China}

\email{nzhang@ouc.edu.cn}

%\tableofcontents

\maketitle

\pagestyle{myheadings} \markboth{\centerline{N. Zhang}}{\centerline{Range decreasing group homomorphisms and holomorphic maps}}

\section{Introduction \label{intro}}

Let $M$ be a positive dimensional $C^{\infty}$ manifold without boundary
and $V$ a finite dimensional $C^{\infty}$ compact connected manifold, possibly with boundary
(all finite dimensional manifolds considered in this paper are second countable).
We fix a smoothness class $\cf=C^{\infty}$, H\"older $\ka$ ($k=0, 1, \cdots$, $0 \le \alpha \le 1$, where $\di C^{k, 0}=C^k$)
or Sobolev $\kp$
($k=1, 2, \cdots$, $1 \le p <\infty$, $kp>\dimv$ or $k = \dimv$, $p=1$).
Then the space $\fvm$ of all $\cf$ maps $V \to M$
is a $C^{\infty}$ Banach/Fr\'echet manifold (see \cite{rh, kr, pa}), to which
we refer as the generalized loop space of $M$ (where $\fsom$ is the loop space of $M$).
Note  that
$$\di C^{\infty}(V, M) \subset \fvm \subset C(V, M).$$
If $M$ is a Lie group $\cg$, then
$\fvg$ is a Lie group under pointwise group operation (see \cite{ps, kr}).
If $M$ is a complex manifold, then $\fvm$ carries a natural complex manifold structure (see \cite{l04, ls}).
The mapping spaces $\fvm$ and operators $\fovn \to \ftwm$
are of fundamental importance in analysis, geometry, mathematical physics and representation theory.

Let $\tovtw$ be the space of maps $\phi: W \to V$ such that
$x \circ \phi \in \ftwr$ for every $x \in \fovr$. It always contains constants.
If $k_1>k_2$, or $k_1=k_2$ and $\alpha_1 \ge \alpha_2$, then all $C^{\infty}$ maps $W \to V$ are in $\tkovktw$.
We write $\di \tfv$ for the space $\tfvfv$.
If $\phi \in \tovtw$, then the pullback operator
\begin{equation*} \label{fpotm}
\phi^{\ast}=\paotm: \fovm \ni x \mapsto x \circ \phi \in  \ftwm
\end{equation*}
is a well defined $C^{\infty}$ map. It is holomorphic when $M$ is a complex manifold
(see Proposition \ref{fp}). If $M$ is a Lie group, it is clear that $\pa$ is a Lie group homomorphism.
Let $\Omega$ be an open subset
  of $\fovm$. We say that a map $f: \Omega \to \ftwm$ is range decreasing if
\begin{equation*} \label{rotation}
  f\left(x\right)\left(W\right) \subset x\left(V\right), \hspace{2mm} x \in \Omega.
\end{equation*}
It is clear that $\di \pa|_{\Omega}$ is range decreasing.
By \cite[Lemma 3.2]{z17}, a range decreasing map $f: \cf(S^1, \bc^n) \to \cf(S^1, \bc^n)$
(where $\cf=C^k$, $k=1, 2, \cdots, \infty$, or $\cf=\kp$) is a pullback operator
if and only if $f$ is continuous complex linear.  This lemma plays a key role
in the proof of \cite[Theorem 1.1]{z17}: if a holomorphic self-map $f$ of $\fsopn$ induces
an isomorphism of the second integral cohomology group of $\fsopn$, then we have the decomposition
\begin{equation*}
  f=\gamma \circ \phi^{\ast}_{\bp^n},
\end{equation*}
where $\phi \in \st(\cf, S^1)$ and $\gamma \in \cf(S^1, PGL(n+1, \bc))$ (which
can be considered as a holomorphic automorphism of $\fsopn$ induced by a family of holomorphic
automorphisms of $\bp^n$ parameterized by $t \in S^1$).
A range decreasing $C^{\infty}$ map is not always (the restriction of)
a pullback operator (see Subsection \ref{pb}). As we shall prove in this paper,
oftentimes the range decreasing condition does imply
that a group homomorphism resp. holomorphic map is (the restriction of) a pullback operator.

The constant maps $V \to M$ form a submanifold of $\fvm$, which can be identified with $M$.
We write $\fvgz$ for the ``trivial'' component of the Lie group $\fvg$ containing the identity;
and if $M$ is connected, we write $\di \fvmz$ for the ``trivial'' component of $\fvm$ containing the constants.
Recall that the evaluation at $v \in V$:
$$E_v=E_{v, \fvm}: \fvm \ni x \mapsto x(v) \in M$$
is a $C^{\infty}$ map. If $M$ is a Lie group resp. complex manifold, then $E_v$ is a Lie group
homomorphism resp. holomorphic map.

In this paper we provide (partial) answers to the following questions:

\begin{question}
  Let $f:$ $\fovgz$ $\to$ $\ftwg$ resp. $f: \fovg \to \ftwg$ be a range decreasing group homomorphism (in the algebraic sense).
  Under what conditions (on $V$, $W$, $\cg$, $\fo$ and $\ft$) must $f$ be (the restriction of) a pullback operator?
\end{question}

\begin{question}
  Let $\Omega$ be a connected open subset
  of $\fovm$ and $f: \Omega \to \ftwm$ a range decreasing holomorphic map.
  Under what conditions (on $V$, $W$, $M$, $\fo$, $\ft$ and $\Omega$) must $f$ be (the restriction of) a pullback operator?
\end{question}

Answers to Questions 1 and 2 are useful to decompose some of the group homomorphisms
and holomorphic maps between generalized loop spaces. As an illustration,
we consider the following

\begin{question}
  Let $\Omega$ be a connected open subset
  of $\fovn$ and $f: \Omega \to \ftwm$ a holomorphic map. How can we tell whether $f$ has a decomposition
  of the form
  \begin{equation} \label{de}
  \di f=\gamma \circ \phi^{\ast}_N|_{\Omega},
  \end{equation}
  where $\phi \in \tovtw$ and $\gamma$ is a holomorphic map from a connected open neighborhood $\tilde{\Omega} \subset \ftwn$
  of $\di \phi^{\ast}_N\left( \Omega \right)$ to $\di \ftwm$ induced by
  a family of holomorphic maps $\gamma_w: N \supset E_w(\tilde{\Omega}) \to M$ with
  $$\di \gamma(y)(w)=\gamma_w \left(y(w) \right), \hspace{2mm} w \in W, \hspace{2mm} y \in \tilde{\Omega} \,?$$
\end{question}

\begin{thm} \label{linear}
Let $f$ be a range decreasing group homomorphism from $\fovgz$ resp. $\fovg$ to $\ftwg$. Suppose that $\dim_{\br} \cg \ge 2$.
Then there exists $\phi \in \tovtw$ such that
$\di f=\pa\left(=\pa_{\cg}\right)$ on any component of $\fovg$ containing an element whose image is nowhere dense in $\cg$.
In particular, $\di f=\pa$ on $\fovgz$; and
   if $\dim_{\br} V < \dim_{\br} \cg$, then $\di f=\pa$ on $\fovg$.
\end{thm}

Theorem \ref{linear}
does not hold if $\dim_{\br} \cg=1$ (see Section \ref{group}).
The automorphism group $Aut(\cg)$ of $\cg$ is a finite dimensional Lie group
if $\pi_0(\cg)$ is finitely generated (see \cite[Theorem 2]{ho}).
In this case, any $\gamma \in \ftwag$ induces a Lie group automorphism of $\ftwg$ resp.
$\ftwgz$ (see Subsection \ref{family}), which is still denoted by $\gamma$.

\begin{cor} \label{sz}
  Let $f:$ $\fovgz$ $\to$ $\ftwg$ be a group homomorphism. Assume
  that $\dim_{\br} \cg \ge 2$ and $\cg$ is connected. Then
  \begin{equation} \label{nonzero}
  \di f(\fovgo) \subset \ftwgo
  \end{equation}
  if and only if there exist $\gf \in \ftwag$ and $\phi \in \tovtw$
  such that
  \begin{equation} \label{ghd}
  \di f=\gf \circ \pa_{\cg}.
  \end{equation}
\end{cor}

It is clear that $\di \gf(w)=E_w \circ f|_{\cg} \in Aut(\cg)$, and the decomposition in (\ref{ghd}) is unique.

Recall that a $C^1$ immersion $x$ from $V$ to a complex manifold $M$ is totally real
if for each $v \in V$, the real subspace $d_v x (T_v V) \subset T_{x(v)} M$ contains no complex subspace
of positive dimension. If there exists a totally real immersion $V \to M$, then
$\dimv \le \dimm$. Any immersion $\di x \in \cf\left(S^1, M \right)$ or  $\di x \in \cf\left([0, 1], M \right)$ (the space of open strings)
 is totally real. For more information about totally real immersions and related
references see \cite[Section 9.1]{fo}.

\begin{thm} \label{main}
  Let $\Omega$ be a connected open subset of $\fovm$ and $f: \Omega \to \ftwm$ a range decreasing holomorphic map.
  Assume that for some $k=1, 2, \cdots$ and $0 \le \alpha \le 1$,
  there is a continuous embedding
  \begin{equation} \label{inclusion}
  \kavm \hspace{1mm} \text{resp.} \hspace{1mm} \civm \subset \fovm
  \end{equation}
  with a dense image.
  Suppose that there exist
  a connected open subset $V_1 \subset V$ and a $\ka/C^{\infty}$ element $y \in \Omega$ such that
  $y|_{V_1}$ is a totally real embedding,
  $\di y(V \setminus V_1) \cap y(V_1)=\emptyset$ and
  $f(y)(W) \cap y(V_1) \not=\emptyset$. Then $f=\pa|_{\Omega}$ for some $\phi \in \tovtw$.
\end{thm}

For all regularities $\fo$ considered in this paper except $\fo=C^{0, \alpha}$, $0< \alpha \le 1$,
there is a continuous embedding with a dense image as in (\ref{inclusion}). Note that \cite[Lemma  3.2]{z17}
is a special case of Theorem \ref{linear} and of Theorem \ref{main}.
The proofs of Theorems \ref{linear} and \ref{main} are different from that of \cite[Lemma  3.2]{z17}.

\begin{cor} \label{cy}
  Let $f$ be as in Theorem \ref{main}. Suppose that the embedding (\ref{inclusion}) (with a dense image)
  and one of the following conditions hold:

  \begin{itemize}
    \item[(a)] there is a  totally real $\ka/C^{\infty}$ embedding $y \in \Omega$;

    \item[(b)]  there is a  totally real $C^{k+\mathrm{sgn}(\alpha)}/C^{\infty}$ immersion $y \in \Omega$;

    \item[(c)] $\di \Omega$ contains a $C^{k+\mathrm{sgn}(\alpha)}/C^{\infty}$ element and
    $\di \dim_{\bc} M \ge \left[\frac{3}{2} \dim_{\br} V \right]$ (the integer part of $\di \frac{3}{2} \dim_{\br} V $).
  \end{itemize}
  Then $f=\pa|_{\Omega}$ for some $\phi \in \tovtw$.
\end{cor}

The map $f$ in Question 3 is completely determined by the family of maps $E_w \circ f: \Omega \to M$, $w \in W$.
If $f$ can be decomposed as in (\ref{de}), then $E_w \circ f=\gamma_w \circ E_{\phi(w)}|_{\Omega}.$
Suppose $N \subset \Omega$. Then we have $\gamma_w=E_w \circ f|_{N}$.
Furthermore, if $E_w \circ f|_{N}$ is not constant for every $w \in W$, then
the decomposition (\ref{de}) is unique.
For any $x \in \fvm$, the pullback $\xtm$ of the tangent bundle of $M$ by $x$
is an $\cf$ bundle (i.e.
the transition functions of the bundle are $\cf$ maps), and
the tangent space $T_x\fvm$ is $\cf(x^{\ast}TM)$, the Banach/Fr\'echet space of
$\cf$ sections of $\xtm$ (see \cite{l04, ls}).
If $x^{\ast}TM$ is trivial, then $T_x\fvm$ can be considered as
$\fvrm$ or $\fvcm$ (where $m=\dim_{\br} M$ or $\dim_{\bc} M$).
As an application of Corollary \ref{sz}, we
obtain the following

\begin{thm} \label{kernel}
  Let $\Omega$ be a connected open subset
  of $\fovn$ and $f: \Omega \to M$ a holomorphic map. Assume that $N, M$ are connected,
  $\dim_{\bc}N=\dim_{\bc} M=n \ge 1$, $\xtn$ is trivial
  for any $x \in \Omega$ and
  \begin{equation} \label{x1x2}
    f(x_1)\not=f(x_2) \hspace{1mm} \text{for} \hspace{2mm} x_1, x_2 \in \Omega \hspace{2mm} \text{with} \hspace{2mm}
     x_1(v) \not= x_2(v), \hspace{2mm} v \in V.
  \end{equation}
  Then there is a fixed $v_0 \in V$ such that
  $$\ker d_x f=\ker E_{v_0, \fovcn}, \hspace{2mm} x \in \Omega$$
  (where $\ker d_x f$ is the kernel of $d_x f$).
  In particular, $f$ is constant on any component of the complex submanifolds $\Omega \cap E_{v_0, \fovn}^{-1}(\zeta)$,
  $\zeta \in N$.
\end{thm}

Submanifolds $E_{v_0, \fovn}^{-1}(\zeta) \subset \fovn$ are of complex codimension $\dim_{\bc} N$ (see \cite[p. 40]{l04}).
Recall the local charts of $\fovn$ as in \cite{l04, ls}.
It is clear that the map $f$ in Theorem \ref{kernel} has the following property: For any $x \in \Omega$,
    there exists a neighborhood $\Omega_x \subset \Omega$ of $x$ such that $f|_{\Omega_x}$ is the composition
    of $E_{v_0}=E_{v_0, \fovn}$ and an injective map $E_{v_0}(\Omega_x) \to M$. For a holomorphic map
    $h: N \to M$, we write $N_{j, h} \subset N$ for the subset of points $x$ such that the complex rank of $d_x h$ is $j$.

\begin{cor} \label{iff}
  Let $O \subset \fovnz$ be a connected open neighborhood of $N$ and $f: O \to M$ a holomorphic map with
  $\dim_{\bc} N=n \le \dim_{\bc} M$.  Then
  $\di f=f|_{N} \circ E_{v_0}|_{O}$ for some $v_0 \in V$ with $N_{n, f|_N} \not=\emptyset$ if and only if
  there exists
  an open subset $\Omega \subset O$ such that $\Omega \cap N \not=\emptyset$, $\di f(\Omega)$ is contained
  in an $n$ dimensional complex submanifold of $M$ and
  (\ref{x1x2}) holds on $\Omega$.
\end{cor}

 By \cite{z03, lz, z10}, there does not exist a $PGL(2, \bc)$-equivariant holomorphic embedding of $\fsopo$ into
 a projectivized Banach/Fr\'echet space $P$ (where the action of $PGL(2, \bc)$ on $P$ is given
 by a monomorphism from $PGL(2, \bc)$ to the group
 of holomorphic automorphisms of $P$).
On the other hand, if $e: M \to \pn$ is a holomorphic embedding, then $e_{\ast}: \fvm \ni x \mapsto e \circ x \in \fvpn$ is a
holomorphic embedding with $e_{\ast} \left(\fvmz \right) \subset \fvpnz$.
If $e$ is $G$-equivariant for a group $G$, then so is $e_{\ast}$.
Thus $\fvpn$ is expected to play a special role
in the theory of generalized loop spaces of complex projective manifolds.

For integers $m \ge n \ge 1$ and $d \ge 0$, let $\hd$ be the space of
holomorphic maps $\bp^n \to \bp^m$ which induce multiplication by $d$ in the second integral cohomology.
The space $\hdo$ consists of constants. The space $\hd$ with $d \ge 1$ can be identified with
a connected open subset of a finite dimensional complex projective space;
and  any element of $\fvhd$ can be considered as a holomorphic map
$\fvpn \to \fvpm$. In particular, $\ho=PGL(n+1, \bc)$ is the group of holomorphic automorphisms
of $\bp^n$, and $\fvpglnp$ acts on $\fvpn$ holomorphically  (see Subsection \ref{family}).
Let $f$ be a holomorphic map from an open neighborhood of $\bp^n \subset \fovpnz$ to $\ftwpm$.
We say that the degree of $f$ is $d$ and write
$$\deg f=d$$
if $E_w \circ f|_{\bp^n} \in \hd$, $w \in W$.
As a consequence of Corollary \ref{iff}, we obtain the following

\begin{thm} \label{degreeone}
Let $f: \fovpnz \to \ftwpm$ be a holomorphic map of degree $d$. Suppose that
one of the following conditions holds:
\begin{enumerate}
  \item[(a)] $d=1$;

  \item[(b)] $d=2$ and $n=m=1$.
\end{enumerate}
Then there exist $\gf \in \ftwhd$ and
  $\phi \in \tovtw$ such that
  \begin{equation} \label{decomposition}
     f=\gf \circ \paotpn.
  \end{equation}
\end{thm}

Thus any range decreasing holomorphic map $f:$ $\fovpnz$ $\to$ $\ftwpn$ must be a pullback operator
($f|_{\bp^n}=\mathrm{id}$, so $\df=1$ and $\gf=\mathrm{id}$).
On the other hand, for any nontrivial component $\Omega_j$ of $\cistpo$,
there exist range decreasing holomorphic maps $\di \Omega_j \to \cistpo$
which are not pullback operators (see Section \ref{holorot}, cf. Corollary \ref{cy}(c)).

Let $f$ be as in (\ref{decomposition}) and $w \in W$. Then $\di (E_w \circ f)(x)$ is completely determined by
the $0$-jet of $x$ at $\phi(w)$.
Given $\phi_1, \cdots, \phi_p \in \tcisotw$ and positive integers $r_1, \cdots, r_p$, we shall construct in Section \ref{pn}
holomorphic maps $g: \cipo \to \ftwpm$ for sufficiently large $m$ such that
$E_w \circ g(x)$ is completely determined by
$r_j$-jets of $x$ at $\phi_j(w)$, $j=1, \cdots, p$, and $E_w \circ g(x)$
can not be completely determined by $0$-jets of $x$ at $\phi_j(w)$, $j=1, \cdots, p$.
In particular, Theorem \ref{degreeone}(b) does not hold if $n=1$ and $m \ge 2$.
These constructions are based on the results in \cite[Section 4]{z03} about holomorphic sections of line bundles
over $\cipo$. For more information about holomorphic sections of line bundles over $\fsopo$ see \cite{z03, z10, ns}.

We denote by $\dfv$ the space of bijections $\phi: V \to V$ such that
$\phi, \phi^{-1} \in \tfv$.
Then for any connected complex manifold $M$, we may consider $\dfv$
as a subgroup of the group $\di \auto$ of holomorphic automorphisms of $\fvmz$.
For more information about $\dfv$ see Subsection \ref{pb}.
It may happen that an element of $\di \fvpglnp$ sends $\fvpnz$ to a different component
of $\fvpn$ (see Section \ref{nm}). Let $\fvpglnptr$ be the open subgroup of $\fvpglnp$
consisting of elements $\gamma$ with $\di \gamma\left(\fvpnz \right)=\fvpnz$.

\begin{cor} \label{auto}
  The group $\autopn$ is the semidirect product $$\fvpglnptr \rtimes \dfv.$$
\end{cor}

The space $\di \dciv$ can be endowed with
a Lie group structure such that the action of $\dciv$ on $\di \civpn$ is $C^{\infty}$ (see \cite[Theorem 11.11 and Remark 11.5]{mi}). Note that
there does not exist a complex structure on $\dciso$ which is compatible with its Lie group
structure (see \cite[Proposition 3.3.2]{ps}).

This paper is organized as follows. Section \ref{bg} contains basic facts about the two classes of maps
in the decomposition (\ref{de}) and about the space $\fvpnz$.

In Sections \ref{tr} and \ref{holorot}, we
give partial answers to Question 2.
For a range decreasing holomorphic map $f: \fovm \supset \Omega \to \ftwm$, define
\begin{equation} \label{dwf}
\wf=\left\{w \in W: \text{there is a} \hspace{1mm} v \in V \hspace{1mm} \text{with} \hspace{1mm} f(x)(w)=x(v) \hspace{1mm} \text{for all} \hspace{1mm} x \in \Omega \right\}.
\end{equation}
It is clear that $\wf$ is a closed subset of $W$.
With the notation in Theorem \ref{main}, it turns out that if $f(y)(w_0) \in y(V_1)$,
  then $w_0$ is in the interior of $\wf$ (Theorem \ref{wf});
  and if $\dimv \le \dimm$, $\wf \not= \emptyset$, then $\wf=W$ (Lemma \ref{wfw}).
Theorem \ref{main} is a consequence of Theorem \ref{wf} and Lemma \ref{wfw}. To show Corollary \ref{cy}, we only
need to find a totally real embedding (when $\dimv<\dimm$) or a totally real immersion
with normal crossings (when $\dimv=\dimm$) in $\Omega$ (Proposition \ref{normal}).

In Section \ref{group}, we try to answer Question 1.
We prove that any group homomorphism in Question 1 is a Lie group homomorphism (Proposition \ref{lgh}).
Let $f$ be as in Theorem \ref{linear}.
To show that $\di f=\pa$ on $\fovgz$,
we pass to the Lie algebra homomorphism $\fh$ induced by $f$ and consider the maps $E_w \circ \fh$, $w \in W$,
so the general case can be reduced
to the special case when $\cg=\br^n$, $\dim W=0$ and $f$ is real linear.
If a real linear map $\di f=(\tilde{f}_1, \cdots, \tilde{f}_n): \fovrn \to \br^n$, where $n \ge 2$, is range decreasing,
then $\di \tilde{f}_1, \cdots, \tilde{f}_n$ must satisfy certain compatibility conditions, which imply that $f$ is a pullback operator (Proposition \ref{rn}).

In Section \ref{nm}, we study the decomposition of holomorphic maps as in Question 3.
For the map $f$ in Theorem \ref{kernel}, we show that its differential $d_x f$
at any $x \in \Omega$ satisfies (\ref{nonzero}) (where $\dim W=0$ and $\cg=\bc^n$).
By Corollary \ref{sz}, $df$ can be considered as a family of maps holomorphically  parameterized by $x \in \Omega$,
and each of them has the decomposition (\ref{ghd}). It turns out that the map $\pa$ in the decomposition
is independent of $x$ (Lemma \ref{tg}). Corollary \ref{iff} is a consequence of Theorem \ref{kernel}.
Some special cases of Corollary \ref{iff} could also be proved by Theorem \ref{main} or by Corollary
\ref{cy}.

In the final Section \ref{pn}, we focus on
holomorphic maps $\di f:$ $\fovpnz$ $\to$ $\ftwpm$. Such a map $\di f$
can be decomposed as in Question 3
if and only if $E_w \circ f$, $w \in W$, have the same kind of decompositions ((\ref{gf}) and Lemma \ref{phi}).
Therefore to prove Theorem \ref{degreeone}, we only need to deal with the special case when $\dim W=0$.
By examining how rational curves in $\fovpnz$ are transformed by $f$,
we show that $f$ has the properties in Corollary \ref{iff} (Proposition \ref{not=}).
Hence we have the decomposition (\ref{decomposition}). At the end of this paper we construct holomorphic maps
$\cipo \to \ftwpm$ which are not of the form (\ref{decomposition}).
Still pullback operators play an important role in these constructions.

We refer to \cite{d, h} for the fundamentals of infinite dimensional
holomorphy.

%%%%%%%%%%%%%%%%%%%%%%%%%%%%%%%%%%%%%%%%%%%%%%%%%%%%%%%%%%%%%%%%%%%%%%%%%%%%%%%%%%%%%%%%%%%%%%%%%%%%%%%%%%%%%%%%%%%

\section{Background \label{bg}}

\subsection{Smooth families of holomorphic maps \label{family}}

Let $Q, N, M$ be finite dimensional $C^{\infty}$ manifolds without boundary
and $\beta: Q \times N \to M$ a $C^{\infty}$ map. Then $\beta$
induces a $C^{\infty}$ map
\begin{eqnarray}
& \bla=\bv: \fvq \times \fvn \to \fvm, \hspace{2mm} \text{where} & \notag \\
& \label{push}
\bla(\gamma, x)(v)=\beta\left(\gamma(v), x(v) \right), \hspace{2mm} v \in V &
\end{eqnarray}
(see \cite[p. 74]{kr}, \cite[p. 91]{rh}). Hence any $\gamma \in \fvq$ determines a $C^{\infty}$ map
$\fvn \to \fvm$, which is still denoted by $\gamma$.
If $N$ and $M$ are complex manifolds
and $\beta(q, \cdot)$ is holomorphic for every $q \in Q$, then $\gamma: \fvn \to \fvm$
is holomorphic
(see \cite[Proposition 2.3]{l04}, where we set $\Phi$ to be
$\di V \times N \ni (v, z) \mapsto \beta(\gamma(v), z) \in M$;
also see \cite[p. 42]{l04} for the case when $\cf=\ka/\kp$).
If $Q, N$ and $M$ are all complex manifolds and $\beta$ is holomorphic, then
$\bla$ is holomorphic. In particular, a $C^{\infty}$ resp. holomorphic
  action of a Lie group $\cg$ on a manifold $M$ induces a $C^{\infty}$ resp. holomorphic
  action of $\fvg$ on $\fvm$. Setting $Q=\{q\}$ (a single point), we see that any $C^{\infty}$ resp. holomorphic map $\di \tilde{\beta}: N \to M$
induces a $C^{\infty}$ resp. holomorphic map $\di \tilde{\beta}_{\ast}: \fvn \ni x \mapsto \tilde{\beta} \circ x \in \fvm$.
Let $\exp: \fg \to \cg$ be the exponential map of a finite dimensional Lie group. Then
$$\di \exp_{\ast}=\exp_{\ast}^{\cf, V}: \fvfg \to \fvg$$ is the exponential map of $\fvg$.
If $U$ is an open neighborhood of $\zero \in \fg$ such that $\exp|_U: U \to \exp(U)$ is a diffeomorphism, then
$$\di \cf(V, U) \ni x \mapsto \exp_{\ast}(x) \in \cf(V, \exp(U))$$ is a diffeomorphism.

Let $\phi \in \tovtw$ and $\gamma \in \fovq$. It is straightforward to verify that
  \begin{equation} \label{comm}
    \paotm \circ \gamma=\paotq(\gamma) \circ \paotn: \fovn \to \ftwm,
  \end{equation}
  where we consider $\gamma$ on the left hand side as  a map $\di \fovn \to \fovm$,
  and consider $\di \paotq(\gamma) \in \ftwq$ as a map $\di \ftwn \to \ftwm$.
  Equation (\ref{comm}) is a generalization of \cite[(2)]{z17}.

For integers $m \ge n \ge 1$ and $d \ge 1$,
let $E=E_d^{n, m}$ be the complex vector space of $(m+1)$-tuples of homogeneous polynomials of degree $d$ on $\bc^{n+1}$.
Then any $h \in \hd$ can be represented by an element of $E$ such that the $m+1$ polynomials have
no common zeros on $\bc^{n+1} \setminus \{\mathbf{0}\}$, and such a representation is unique up to an overall multiplicative constant.
Thus $\hd$ can be considered as an open subset of the projective space $P(E)$,
and the evaluation
\begin{equation*}
  \hd \times \pn \ni (h, \zeta) \mapsto h(\zeta) \in \bp^m
\end{equation*}
is holomorphic.
So the induced map
$$\fvhd \times \fvpn \to \fvpm$$
as in (\ref{push}) is holomorphic.
The space $\hd$ is connected,
and it is simply connected if $n<m$. For more information about the topology of $\hd$ see \cite{mo, ya, fe}.

Suppose $\gamma: V \to \hd$ is a map such that $\gamma(\cdot)(\zeta): V \to \bp^m$ is an $\cf$ map
for any fixed $\zeta \in \pn$.
By induction on $n$, one can verify that for any $\zeta_0 \in \bp^n$ and any $v_0 \in V$,
there exist a neighborhood $U_0$ of $\zeta_0$ and a neighborhood $V_0$ of $v_0$
on which $\gamma(v)(\zeta)$ can be represented (under affine coordinates of $\bp^n$ and of $\bp^m$)
by $m$-tuples of rational functions of $\zeta \in U_0$
whose coefficients are $\cf$ functions of $v \in V_0$. So
$\di \gamma \in \fvhd$.
Let $x \in \fvpn$, and let
\begin{equation} \label{kappa}
\kappa=\kappa_{\gamma, x}: x^{\ast}T\bp^n  \to \gamma(x)^{\ast}T\bp^m
\end{equation}
be the morphism of vector bundles covering the identity map on the base space $V$ such that the restriction
of $\kappa$ to the fiber over $v \in V$ is $\di d_{x(v)} \gamma(v)$. Note that
$\kappa$ is an $\cf$ morphism of vector bundles (i.e. under local trivializations of $x^{\ast}T\bp^n$
and of $\gamma(x)^{\ast}T\bp^m$, $\kappa$
can be represented by $\cf$ maps from open subsets of $V$ to the space of $m \times n$ complex matrices).
It is clear that $\di E_{v, \fvpm} \circ \gamma=\gamma(v) \circ E_{v, \fvpn}.$ So we have
$E_{v, \fgxtpm} \circ d_x \gamma=d_{x(v)}\gamma(v) \circ E_{v, \fxtpn},$
where $\di  \fxtpn$ (resp. $\di \fgxtpm$) is the tangent space
of $\fvpn$ (resp. $\fvpm$) at $x$ (resp. $\gamma(x)$). Thus
\begin{equation} \label{dx}
\di d_x \gamma(\tx)=\kappa \circ \tx, \hspace{2mm} \tx \in \fxtpn.
\end{equation}

\subsection{Basic properties of pullback operators \label{pb}}

Let $\phi \in \tovtw$. For sufficiently large $n$, there is a $C^{\infty}$ embedding $x_0: V \to \br^n$.
So $\phi=x_0^{-1} \circ \left(x_0 \circ \phi \right)$ is an $\ft$ map.

\begin{prop} \label{fp}
Suppose $M$ is an $m$ dimensional $C^{\infty}$ manifold without boundary, where $m=1, 2, \cdots$,
$\phi: W \to V$ is a map and $\Omega \subset \fovm$ is a nonempty open subset. Then $\phi \in \tovtw$ if
and only if $x \circ \phi \in \ftwm$ for every $x \in \Omega$. In this case,
the map $\pa_M$ is $C^{\infty}$. Furthermore, if $M$ is a complex manifold,
then $\pa_M$ is holomorphic.
\end{prop}

\begin{proof}
  To show the ``only if'' direction, let $(U, \Phi)$ be a coordinate chart of $M$ and $x \in \fo(V, U) \subset \fovm$. Then
  \begin{equation*}
    x \circ \phi=\iptw \circ \paotrm \circ \pov(x) \in \ft(W, U) \subset \ftwm.
  \end{equation*}
  With cut-off functions one can show that
  $x \circ \phi \in \ftwm$ for every $x \in \fovm$ (see the proof of \cite[Proposition 2.1]{z17}).

Regarding the ``if'' direction, recall the local chart $(\su, \varphi_y)$ of $\fovm$ as in
\cite[Subsection 1.1]{ls} (if
$M$ is a $C^{\infty}$ manifold instead of a complex manifold, then similar
 arguments yield a $C^{\infty}$ local chart rather than a holomorphic one), which maps an open neighborhood $\su$ of
$\di y \in \fovm$ to an open neighborhood of $\varphi_y(y)=\zero \in \fo(y^{\ast}TM)$.
Choose $\su \subset \Omega$. We write $\di \pa_M|_{\su}$ for the map $\di \su \ni x \mapsto x \circ \phi \in \ftwm$.
It is straightforward to verify that
$\di \varphi_{y \circ \phi} \circ \pa_M|_{\su} \circ \varphi_y^{-1}$ is given by
\begin{equation*} \label{phia}
    \pa_{y^{\ast}TM}|_{\varphi_y(\su)}:   \varphi_y(\su) \ni \sigma \mapsto \sigma \circ \phi \in \ft(\phi^{\ast}y^{\ast}TM).
\end{equation*}
The linear operator $\di \pa_{y^{\ast}TM}$ is actually
well defined on all of $\fo(y^{\ast}TM)$. Let $\tx \in \fovr$ and $w \in W$. Take
$\di \sigma \in \fo(y^{\ast}TM)$ with $\sigma \circ \phi(w)\not=\zero$. Then
$\di \tx \circ \phi=\pa_{y^{\ast}TM}(\tx \sigma)/\pa_{y^{\ast}TM}(\sigma)$ is an $\ft$ function in a neighborhood of $w$. Therefore $\phi \in \tovtw$.

If a sequence $\{\sigma_n\}$ in $\fo(y^{\ast}TM)$ converges to $\zero$, then the sequence $\{\sigma_n \circ \phi\}$ in $C(\phi^{\ast}y^{\ast}TM)$
converges to $\zero$.
By Closed Graph Theorem, $ \pa_{y^{\ast}TM}$ is continuous. Thus $\pa_M$ is $C^{\infty}$. If $M$
is a complex manifold, then $ \pa_{y^{\ast}TM}$ is complex linear. So $\pa_M$ is holomorphic. \end{proof}

Let $\Omega$ be as in Proposition \ref{fp} and $f: \Omega \to \ftwm$ a map.
Then $f=\pa_M|_{\Omega}$ for some $\phi \in \tovtw$ if and only if for any $w \in W$,
there exists $v(w) \in V$ with
\begin{equation} \label{evw}
E_{w, \ftwm} \circ f=E_{v(w), \fovm}|_{\Omega}.
\end{equation}

Suppose $\phi: V \to V$ is a bijection. If $\phi \in \dfv$, then $\phi, \phi^{-1} \in \cf(V, V)$.
The converse is also true if $\cf$ is one of the following regularities: $C$, $\ka$ with $k+\alpha \ge 1$,
$\di W^{k, p}$ with $k=2, 3, \cdots$, $p>\dimv$ and $\di W^{k, 1}$ with $\dimv=1$, $k=1, 2, \cdots$ (see \cite[Lemma 2.3]{bhs}).
If $1<p <\infty$ and $\dimv=1$,
then $\mathscr{D}\left(W^{1, p}, V\right)$
is the space of bi-Lipschitz maps (see \cite[p. 710]{z17}).

A range decreasing $C^{\infty}$ map is not always a pullback operator. For example, the map
$$C^{\infty}(S^2, \bc^n) \times C^{\infty}(S^2, S^2) \ni (x, y) \mapsto x \circ y \in C^{\infty}(S^2, \bc^n)$$
is $C^{\infty}$ (see \cite[p. 91]{rh}). Let $\tau: C^{\infty}(S^2, \bc^n) \to \bc$
be a non-zero continuous complex linear functional and let $T_{\zeta} \in C^{\infty}(S^2, S^2)$ be the
translation $\bp^1 \ni z \mapsto z+\zeta \in \bp^1$, where $\zeta \in \bc$. Then
$$f_0: C^{\infty}(S^2, \bc^n) \ni x \mapsto x \circ T_{\tau(x)} \in C^{\infty}(S^2, \bc^n)$$
is a range decreasing $C^{\infty}$ map. It is not complex linear (the property $f_0(\lambda x)=\lambda f_0(x)$, where $\lambda \in \bc \setminus \{0\}$,
implies that $x \circ T_{\lambda \tau(x)}=x \circ T_{\tau(x)}$).
So it is not a pullback operator.

\subsection{The space $\fvpn$ \label{geo}}

The sheaf of germs of complex valued $\cf=C^{\infty}/\ka/\kp$ functions over $V$ is fine.
So the isomorphism class of an $\cf$ complex  line bundle $E \to V$ is completely determined by its first Chern class
$c_1(E) \in \htv$. Let $\xi$ be the universal line bundle over $\pn$ and $x \in \fvpn$.
Then $c_1(x^{\ast} \xi)$ is constant on any component of $\fvpn$.
Let $\cu=\cu(\fvpn)$ be the set $\{\mu(\bp^1)\}$ of curves  in
$\fvpn$, where $\mu$ ranges over all holomorphic embeddings $\bp^1
\to \fvpn$ such that $E_v \circ \mu \in \hyn$,
$v \in V$.

\begin{prop} \label{2points}
  Let $x_1, x_2, x_3 \in \fvpn$ be three different elements.
  \begin{itemize}
    \item[(a)] If $\di c_1(x_1^{\ast} \xi)=c_1(x_2^{\ast} \xi)=0$, then there is
    a curve in $\cu$ through both $x_1$ and $x_2$ if and only if $x_1(v)
    \not= x_2(v)$ for all $v \in V$.

    \item[(b)] If $n=1$ and $x_i(v) \not= x_j(v)$ for all $v \in V$, $i, j=1, 2, 3$, $i \not=j$,
    then there is a curve in $\cu$ through $x_1, x_2$ and $x_3$.
  \end{itemize}
\end{prop}
\begin{proof}
  (a) Let $\pi: \bc^{n+1} \setminus \{\zero\} \to \pn$ be the natural projection.
  Then $\di c_1(x^{\ast}_j \xi)=0$, $j=1, 2$, if and only if there exist $\cf$ maps $\tx_j: V \to \bc^{n+1} \setminus \{\zero\}$ with $\pi \circ \tx_j=x_j$.
  The conclusion of (a)
  follows from similar arguments as in the proof of \cite[Proposition 3.3]{z17}.

    (b) Let
    $\gamma$ be the element of $\fvpglt$
    such that for any $v \in V$, $\gamma(v)$ is the unique element of $PGL(2, \bc)$ which
    maps $0, \infty, 1 \in \bp^1$ to $x_1(v)$, $x_2(v)$ and $x_3(v)$ respectively.
    Consider $\gamma$ as a holomorphic automorphism of $\fvpy$. Then $\gamma(\bp^1)$ (where $\bp^1 \subset \fvpy$)
    is a curve in $\fc$ through $x_1, x_2$ and $x_3$.
\end{proof}

If $\dimv \le 2n$, then the map $x \to c_1(x^{\ast} \xi)$ induces a bijection
$$\pi_0 \left(\fvpn \right) \to \htv$$
(e.g. see \cite[Lemma 1.1]{fe} and \cite[Theorem 13.14]{pa}).
In particular, $\fvpnz$ consists of elements $x$ with $c_1(x^{\ast} \xi)=0$.
If $\di \dimv>2n$, then in general the above bijection does not hold. For
example, $H^2(S^3, \bz)=0$, and it follows
from the relation between pointed and free homotopy classes
that $\pi_0 \left(\cf(S^3, \bp^1) \right) \simeq  \pi_3(S^2)/\pi_1(S^2) \simeq \bz$.

%%%%%%%%%%%%%%%%%%%%%%%%%%%%%%%%%%%%%%%%%%%%%%%%%%%%%%%%%%%%%%%%%%%%%%%%%%%%%%%%%%%%%%%%%%%%%%%%%%%%%%%%%%%%%%%%%%%%%%%%%%%%%%

\section{Total reality and the interior of $\wf$ \label{tr}}

The following is the main result of this section.

\begin{thm} \label{wf}
  Let $f$, $\fovm$, $y$, $V_1$ be as in Theorem \ref{main} and $\wf$ be as in (\ref{dwf}).
  If $f(y)(w_0) \in y(V_1)$, where $w_0 \in W$,
  then $w_0 \in \wf^o$ (the interior of $\wf$).
\end{thm}

\begin{prop} \label{continuity}
  Let $M$ be a positive
  dimensional topological manifold, $V$ a compact topological manifold, possibly with boundary,
  $\sq$ a Hausdorff topological space,
  $\sq \ni \lambda \mapsto \xl \in C(V, M)$ (with the compact-open topology) a continuous map, $\tv \subset V$ an
  open subset and $\di \sv: \sq \to \tv$ a map such that $\di \sq \ni \lambda
  \mapsto \xl\left(\sv(\lambda)\right) \in M$ is continuous. If for any $\lambda \in \sq$, the map $\xl$ is
  injective on the closure $\ctv$ of $\tv$, then $\di \sv$ is continuous.
\end{prop}
\begin{proof} Fix a metric on $M$. For any $\lambda_0 \in \sq$ and any open neighborhood $V_1 \subsetneqq \tv$ of $\sv(\lambda_0)$,
we have
$$\di \delta=\mathrm{dist} \left(\xlz\left(\sv(\lambda_0)\right), \xlz(\ctv \setminus V_1) \right)>0.$$
  Since the maps $\sq \ni \lambda \mapsto \xl \in C(V, M)$ and $\di \sq \ni \lambda
  \mapsto \xl\left(\sv(\lambda)\right) \in M$ are continuous, there is
  an open neighborhood $\sq_1$ of $\lambda_0$ such that
  \begin{eqnarray*}
  & \di \mathrm{dist} \left(\xlz\left(\sv(\lambda_0)\right), \xl(\ctv \setminus V_1)\right) >\delta/2 \hspace{2mm} \text{and} & \\
  & \mathrm{dist} \big(\xlz\left(\sv(\lambda_0)\right), \xl\left(\sv(\lambda)\right) \big)<\delta/2 &
  \end{eqnarray*}
  for all $\lambda \in \sq_1$.
  So we have $\di \sv(\sq_1) \subset V_1$, and $\sv$ is continuous at $\lambda_0$.
\end{proof}

\begin{prop} \label{localconstant}
  Let $\di \Delta_r=\left\{\lambda \in \bc: |\lambda|<r \right\}$ ($r>0$),
  $\tv$ a $C^{1}$ manifold, possibly with boundary,
    $\tx, \ty \in \cotvcm$ and
  $\sv: \Delta_r \to \tv$ a continuous map such that the function $h: \Delta_r  \to \bc^m$ defined by
  \begin{equation*} \label{h}
    h(\lambda)=\pxly\left(\sv(\lambda)\right)
  \end{equation*}
  is holomorphic.
  If $\txly: \tv \to \bc^m$ is a $C^1$ embedding
  for every $\lambda \in \Delta_r$, then $\sv$ is real differentiable on $\di \Delta_r$. Furthermore, if
  $\di \txly$ is totally real for every $\lambda \in \Delta_r$, then $\sv$ is constant.
\end{prop}
\begin{proof} For any $\di \lambda \in \Delta_r$,
  \begin{eqnarray*}
    h'(\lambda) & = &  \lim_{\bc \ni \zeta \to 0} \left(\frac{\pxly\left(\sv(\lambda+\zeta)\right)-\pxly\left(\sv(\lambda)\right)}{\zeta}
    +\tx\left(\sv(\lambda +\zeta)\right) \right) \notag\\
    &=& \lim_{\bc \ni \zeta \to 0} \left(\frac{\rho(\zeta)
    -\rho(0)}{\zeta}\right)
    +\tx\left(\sv(\lambda)\right), \label{rho}
  \end{eqnarray*}
    where
  $\rho(\zeta)=\pxly\left(\sv(\lambda+\zeta)\right)$
  is a map from a neighborhood of $0 \in \bc$ to $\bc^m$. It is clear that $\rho$ is complex differentiable at $\zeta=0$. So
    \begin{equation*} \label{tl}
    \sv(\lambda+\zeta)=\left(\txly \right)^{-1}\left(\rho(\zeta)\right)
    \end{equation*}
    is real differentiable at $\zeta=0$.

    The image of $d_0\rho$ is a $J$-invariant subspace of $\di T_{\rho(0)}\pxly(\tv)$, where $J$ is the almost complex structure
    of $\bc^m$. If $\di \txly$ is totally real, then $d_0\rho=\zero$, which implies that $d_{\lambda} \sv=\zero$.
    Thus $\sv$ is constant.
    \end{proof}

\begin{lemma} \label{coordinates}
  Let $M$ be an $m$ dimensional complex manifold without boundary, where
  $m=1, 2, \cdots$, $y \in \fvm$ and $v_0 \in V$. Then there exist a coordinate chart $\di (U, \Psi)$
  of $M$ with $y(v_0) \in U$ and a coordinate
  chart $(\su, \varphi_y)$ of $\fvm$ with $y \in \su$
  such that

  \begin{itemize}

  \item[(i)] $\di \varphi_y(\su) \subset \fytm$,  $\di \varphi_y(y)=\zero$,
  $\varphi_y(x)(v) \in T_{y(v)} M$ for all $x \in \su$ and $v \in V$;

  \item[(ii)] $\Psi(U)$ is a convex open subset of $\bc^m$;

  \item[(iii)] if both $x(v)$ and $y(v)$ are contained in $U$, then
  \begin{equation} \label{substraction}
d\Psi \left(\varphi_y(x)(v)\right)=\left(\Psi\left(x(v)\right)-\Psi\left(y(v)\right),
\Psi\left(y(v)\right)\right) \in \bc^m \times \Psi(U),
  \end{equation}
  where we
  identify the tangent bundle of $\Psi(U)$ with $\bc^m \times \Psi(U)$.
  \end{itemize}
\end{lemma}
\begin{proof} Let $\Psi$ be a biholomorphic map from an open neighborhood of $y(v_0)$ to an open subset of $\bc^m$
with $\di \Psi\left(y(v_0) \right)=\zero$. For a sufficiently small ball $B_{r}(\zero) \subset \bc^m$, define
 $U_0=\Psi^{-1} \left(B_{r}(\zero) \right)$ and $U=\Psi^{-1} \left(B_{r/2}(\zero) \right)$. Take an open covering $\{U_i\}_{i \in \Lambda}$ of $M$
such that every $U_i$ is biholomorphic to a convex open subset of $\bc^m$, $U_0 \in \{U_i\}_{i \in \Lambda}$
and $U_i \cap U=\emptyset$ if $U_i \not=U_0$. Define the diffeomorphism $F$ as in the proof of \cite[Lemma 1.1]{ls}
with the above open covering of $M$.  Shrinking $U$ if necessary, we have
\begin{equation*}
  F(z, w)=(d \Psi)^{-1} \left(\Psi(z)-\Psi(w), \Psi(w) \right), \hspace{2mm} z, w \in U.
\end{equation*}
Construct the coordinate chart $(\su, \varphi_y)$ of $\fvm$ as on page 486 of \cite{ls} with the above diffeomorphism $F$.
It is clear that we have (i), (ii) and (iii).
\end{proof}

\noindent {\it Proof of Theorem \ref{wf}.}
  By (\ref{inclusion}), we may assume that $\fo=\ka/C^{\infty}$.
  It is clear that there is a unique $v_0 \in V_1$ with $f(y)(w_0)=y(v_0)$.
  Let $\di (U, \Psi)$ and $(\su, \varphi_y)$ be as in Lemma \ref{coordinates}.
  Choose a connected open neighborhood $\tv$ of $v_0$ such that
  $\overline{\tv} \subset V_1$, $\overline{\tv}$ is a $C^{\infty}$ submanifold of $V$, possibly with boundary,
  and $\di y\left(\overline{\tv} \right) \subset U$.
  We claim that there exist a connected open neighborhood $\su_1 \subset \Omega \cap \su$ of $y$
  and a connected open neighborhood $\tw \subset W$ of $w_0$ such that for every $x \in \su_1$, we have
  \begin{equation} \label{px}
    \begin{array}{rl}
      (a) & x|_{\overline{\tv}} \,\, \mathrm{is \,\,a \,\, totally \,\,real \,\, embedding}; \\
      (b) & x(\ctv) \subset U; \\
      (c) & f(x)(\ctw) \subset x(\tv). \\
    \end{array} %\left\{ \right.
  \end{equation}
  The pullback operator $i^{\ast}: \fovm
  \to \footvm$ induced by the inclusion $i: \overline{\tv} \to V$ is holomorphic. Note that $\di i^{\ast} y=y|_{\overline{\tv}}$ is a
  totally real embedding and the subset of totally real embeddings
  is open in $\footvm$. So there is a sufficiently small connected open
  neighborhood $\su_1 \subset \Omega \cap \su$ of $y$ such that (a) and (b) hold for every $x \in \su_1$.
  Fix a metric on $M$.
  Since $f(y)(w_0)=y(v_0) \not\in y(V \setminus \tv)$,
  there is a connected open neighborhood $\tw \subset W$ of $w_0$ such that  $\di \mathrm{dist} \left(f(y)(\ctw), y(V \setminus \tv) \right)>0$.
  At the price of shrinking $\su_1$, we can arrange that $\varphi_y(\su_1) \subset \foytm$ is an open ball
  centered at $\zero$ and
  $\di \mathrm{dist} \left(f(x)(\ctw), x(V \setminus \tv) \right)>0$ for every $x \in \su_1$.
  By the range decreasing condition, we obtain (c).

   Next we show that $\di f(x)(w_0)=x(v_0)$ for all $x \in \su_1$, which implies that
  $\di f(x)(w_0)=x(v_0)$ for all $x \in \Omega$,
  hence $w_0 \in \wf$. Note that any $w \in \tw$ satisfies the same kind of conditions
  as $w_0$. So we have $\tw \subset \wf$ and $w_0 \in \wf^o$.
  For any $x \in \su_1 \setminus \{y\}$, let
  \begin{equation*} \label{xl}
    \xl=\varphi^{-1}_y \left(\lambda \varphi_y(x) \right) \in \su_1, \,\, \lambda \in \overline{\Delta},
  \end{equation*}
  where $\di \Delta=\left\{\lambda \in \bc: |\lambda|<1 \right\}$. Then $x_0=y$. By (a) and (c) of (\ref{px}), there is a unique
   $\sv: \overline{\Delta} \to \tv$ with
  \begin{equation} \label{holo}
    f\left(\xl \right)(w_0)=\xl\left(\sv(\lambda)\right): \overline{\Delta} \to M
  \end{equation}
  (where $\sv(0)=v_0$). The above map is holomorphic on $\Delta$
  and continuous on $\di \overline{\Delta}$.
  By Proposition \ref{continuity}, $\sv$ is
  continuous on $\overline{\Delta}$. Let $P_1$ be the projection
  from $\bc^m \times \Psi(U)$ to $\bc^m$. In view of (b) of (\ref{px}) and  (\ref{substraction}),
  \begin{eqnarray*} \label{tyv}
    & \ty(v)=\Psi\left(y(v)\right): \tv \to \bc^m, & \\
     \label{txv}
    & \tx(v)=\left(P_1 \circ d\Psi\right) \left(\varphi_y(x)(v)\right)=\Psi\left(x(v)\right)-\ty(v): \tv \to \bc^m &
  \end{eqnarray*}
  are well defined $\fo$ maps, and
  \begin{eqnarray*}
    \Psi\left(\xl(v)\right)-\ty(v) &=& \left(P_1 \circ d\Psi\right) \left(\varphi_y(\xl)(v)\right)
    =\left(P_1 \circ d\Psi\right) \left(\lambda \varphi_y(x)(v)\right) \\
    &=& \lambda \left(P_1 \circ d\Psi\right) \left(\varphi_y(x)(v)\right)=\lambda \tx(v).
  \end{eqnarray*}
 By (a) of (\ref{px}),
    $\di \ty+\lambda \tx=\Psi \circ \xl: \tv \to \bc^m$
  is a totally real embedding for every $\lambda \in \Delta$.
   By (\ref{holo}), the map
  \begin{equation*}
    \Delta \ni \lambda \mapsto \pxly\left(\sv(\lambda)\right)=\Psi \left( f\left(\xl \right)(w_0) \right) \in \bc^m
  \end{equation*}
  is holomorphic. It follows from Proposition \ref{localconstant} that $\sv$ is constant.
  So $\sv\equiv v_0$ on $\di \overline{\Delta}$ and $\di f(x)(w_0)=x(v_0)$.
\qed

%%%%%%%%%%%%%%%%%%%%%%%%%%%%%%%%%%%%%%%%%%%%%%%%%%%%%%%%%%%%%%%%%%%%%%%%%%%%%%%%%%%%%%%%%%%%%%%%%%%%%%%%%%%%%%%%%%%%%%%%%%%%%%%%%%%%%%%%%%%%%%%%%%%%%

\section{Range decreasing holomorphic maps \label{holorot}}

In this section, we provide partial answers to Question 2. In particular, we prove Theorem \ref{main} and Corollary \ref{cy}.
Let $M$ be a positive dimensional $C^{\infty}$ resp. complex manifold without boundary,
$v_0 \in V$ and $J^k(v_0, M)$, $k=1, 2, \cdots$, the submanifold of the $k$-jet space $J^k(V, M)$
consisting of $k$-jets at $v_0$. It is straightforward to verify that $J^k(v_0, M)$
is a $C^{\infty}$ resp. complex manifold, the target map
$$\tau_k: J^k(v_0, M) \ni j_{v_0}^k x \mapsto  x(v_0) \in M$$
gives rise to a $C^{\infty}$ resp. holomorphic fiber bundle
  and the map $$j^k_{v_0}: \kavm \hspace{1mm} \text{resp.} \hspace{1mm} \civm \ni x \mapsto j^k_{v_0} x \in J^k(v_0, M)$$ is $C^{\infty}$ resp. holomorphic.

\begin{prop} \label{y0}
  Let $M$ be a $C^{\infty}$ manifold without boundary with $\dim_{\br} M$ $=p=2$, $3$, $\cdots$,
  $D \subset \kavm \hspace{1mm} \text{resp.} \hspace{1mm} \civm$ a nonempty open subset, where $k=1, 2, \cdots$,
  and $v_0 \in V$. If $\dimv=l <p$, then
  there exist $y \in D$ and an open neighborhood $V_1$ of $v_0$ such that
  $y|_{V_1}$ is an embedding and
  $\di y(V \setminus V_1) \cap y(V_1)=\emptyset$.
\end{prop}
\begin{proof}
   Let $X_l \subset J^1(v_0, M)$ be the subset of 1-jets $j^1_{v_0} x$ such that
   the rank of $d_{v_0} x$ is $l$.
   For any $\di y_1 \in \kavm$ resp. $\civm$, it follows from Lemma \ref{coordinates} that there is a coordinate
  chart $(\su, \varphi_{y_1})$ of $\kavm$ resp. $\civm$ with $y_1 \in \su$ such that $\di \varphi_{y_1}(\su \setminus \left(j^1_{v_0} \right)^{-1}(X_l))$
  is a real analytic subset of $\varphi_{y_1}(\su)$.
   So $\left(j^1_{v_0} \right)^{-1}(X_l)$
   is an open and dense subset of $\kavm$ resp. $\civm$.

   Take $y_0 \in D \cap \left(j^1_{v_0} \right)^{-1}(X_l)$, a $C^{\infty}$ coordinate chart $(U, \Phi)$
   of $M$ and an open connected neighborhood $\tv$ of $v_0$ such that $\overline{\tv}$ is a $C^{\infty}$ submanifold
   of $V$, possibly with boundary,
   $y_0|_{\overline{\tv}}$ is an embedding and $y_0 \left( \overline{\tv}\right) \subset U$.
   Let $\xi_0 \in S^{p-1}$ and $\chi \in C^{\infty}(V)$ be such that
   $\chi \equiv 1$ on a neighborhood of $v_0$ and $\supp \chi \subset \tv$.
   For a sufficiently small $\ve>0$ and a sufficiently small neighborhood $O \subset S^{p-1}$ of $\xi_0$,
   define $\yxt \in D$ by
   $$\yxt(v)=\left\{ \begin{array}{ll}
     \Phi^{-1} \big( \Phi\left(y_0(v)\right)+t \chi(v) \xi \big), & v \in \tv, \\
     y_0(v), & v \in V \setminus \tv,
   \end{array}\right.$$
   where $\xi \in O$ and $t \in (-\ve, \ve) \subset \br$.  Shrinking $(-\ve, \ve)$ and $O$ if necessary, we may assume that
   $\yxt|_{\overline{\tv}}$ is an embedding for every $ t \in (-\ve, \ve)$ and every $\xi \in O$.
   Note that the interior of the set
   $$\di \left\{\Phi\big(\yxt(v_0)\big) \in \br^{p}: t \in (-\ve, \ve), \xi \in O \right\}$$
   is non-empty and $\Phi \left(U \cap y_0(V \setminus \tv) \right) \subset \br^p$ is nowhere dense.
   Thus there exists $y=\yxoto$ with $y(v_0) \not\in y(V \setminus \tv)$.
   Take an open neighborhood $V_1 \subset \tv$ of $v_0$ such that $y(V_1) \cap y(V \setminus \tv)=\emptyset$. Then
   $y|_{V_1}$ is an embedding and $y(V_1) \cap y(V \setminus V_1)=\emptyset$.
\end{proof}

\begin{lemma} \label{wfw}
  Let $f$ be as in Theorem \ref{main}. Assume that we have the embedding (\ref{inclusion}) (with a dense image).
  If $\dimv=l \le m=\dimm$ and $\wf \not= \emptyset$, then $\wf=W$.
\end{lemma}
\begin{proof}
  Recall that $\wf$ is a closed subset of the connected manifold $W$. We shall prove that
  any $w_0 \in \wf$ is an interior point. Suppose that $f(x)(w_0)=x(v_0)$ for all $x \in \Omega$.
  By Theorem \ref{wf}, we only need to show
  that  there exist
  a connected open neighborhood $V_1$ of $v_0$ and a $\ka/C^{\infty}$ element $y \in \Omega$ such that
  $y|_{V_1}$ is a totally real embedding and
  $\di y(V \setminus V_1) \cap y(V_1)=\emptyset$.

  Let $X \subset J^1(v_0, M)$ be the subset of 1-jets $j^1_{v_0} x$ such that
   $d_{v_0} x$ is injective and totally real. Then $J^1(v_0, M) \setminus X$ is a real
   analytic subset (a real linear map $\br^l \to \bc^m \simeq \br^{2m}$ represented by a
   $2m \times l$ matrix $B$ is injective and totally real if and only if the rank of the $2m \times 2l$
   matrix $(B, JB)$ is $2l$, where $J$ is the standard almost complex structure on $\br^{2m}$).
   Therefore the subset
   $$\left(j^1_{v_0} \right)^{-1}(X) \subset \kavm \hspace{1mm} \text{resp.} \hspace{1mm} \civm$$
   is open and dense.
   By Proposition \ref{y0}, there exist $y \in \left(j^1_{v_0} \right)^{-1}(X) \cap \Omega$
   and an open neighborhood $V'_1$ of $v_0$ such that
  $y|_{V'_1}$ is an embedding and
  $\di y(V \setminus V'_1) \cap y(V'_1)=\emptyset$. Choose a connected open neighborhood $V_1 \subset V_1'$
  of $v_0$ such that $y|_{V_1}$ is totally real, then
  $\di y(V \setminus V_1) \cap y(V_1)=\emptyset$.
\end{proof}

Theorem \ref{main} immediately follows from Theorem \ref{wf}, Lemma \ref{wfw} and Proposition \ref{fp}.

\begin{prop} \label{normal}
  Let $f$ be as in Theorem \ref{main}. Assume that we have the embedding (\ref{inclusion}) (with a dense image).
  If $\dimv=\dimm$ and there exists a totally real $\ka/C^{\infty}$ immersion $y \in \Omega$
  with normal crossings (i.e. the only self-intersections of $y$ are transversal double points),
  then $f=\pa|_{\Omega}$ for some $\phi \in \tovtw$.
\end{prop}
\begin{proof}
  We only need to consider the special case when $\di \fo=\ka/C^{\infty}$.
  Assume for contradiction that $f$ is not a pullback operator. It follows from Theorem \ref{main} that
  for any totally real immersion $x \in \Omega$ and any $w \in W$, $(E_w \circ f)(x)$
  is a self-intersection point of $x$. In particular, $f(y) \in M \subset \ftwm$ is one of
  the double points of $y$.
  Suppose that
  $$f(y)=y(v_0)=y(v_1), \hspace{1mm} v_0, v_1 \in V, \hspace{1mm} v_0 \not= v_1.$$
  Let $V_0 \subset V$ be an open connected neighborhood of $v_0$
  such that $\di \overline{V_0}$ is a $C^{\infty}$ submanifold of $V$, possibly with boundary,
  and $\di y|_{\overline{V_0}}$ is an embedding, $\tx \in \kaytm$ resp. $\ciytm$ with
  $\supp \tx \subset V_0$ and
    \begin{equation} \label{rank}
    \tx(v_0) \in T_{y(v_0)}M \setminus d_{v_0}y \left(T_{v_0} V \right)
    \end{equation}
  (where we consider $\tx(v)$ as an element of $T_{y(v)} M$) and let
  $(\su, \varphi_y)$ be a coordinate chart of $\kavm$ resp. $\civm$ as in Lemma \ref{coordinates} with
  $\su \subset \Omega$.
  Choose
  $\delta>0$ such that for any
  $\lambda \in \Delta_{\delta}=\{z \in \bc: |z|<\delta\}$,
  we have $\lambda \tx \in \varphi_y(\su)$, $\xl=\varphi_y^{-1}(\lambda \tx)$
  is a totally real immersion and $\di \xl|_{\overline{V_0}}$ is an embedding.
   Note that $(E_w \circ f)(\xl)$ is a self-intersection point
  of $\xl$. Thus it is contained in the totally real immersed submanifold
  $\xl(V \setminus \overline{V_0})=y(V \setminus \overline{V_0}) \subset M$,
  which implies that the holomorphic map $\Delta_{\delta} \ni \lambda \mapsto (E_w \circ f)(\xl) \in M$ is constant.
  So $f(\xl)=y(v_0)$ for every $\lambda \in \Delta_{\delta}$.

  On the other hand, it follows from (\ref{rank}) and (\ref{substraction}) that the rank of the map
  $$V \times \left(-\delta, \delta \right) \ni (v, \lambda) \mapsto \xl(v) \in M$$
  at the point $(v_0, 0)$ is $\dimv+1$. So there exist a positive $\delta_1 \le \delta$ and
  an open neighborhood $\tv_0 \subset V_0$ of $v_0$ with $y(v_0) \not\in \xl(\tv_0)$ for every $\lambda \in \left(0, \delta_1\right)$.
  As $y(v_0) \not\in y(\overline{V_0} \setminus \tv_0)$, there is a positive $\delta_2 \le \delta_1$
      such that $y(v_0) \not\in \xl(\overline{V_0} \setminus \tv_0)$ for every $\lambda \in \left(0, \delta_2\right)$.
       Therefore $y(v_0) \not\in \xl(\overline{V_0})$. Note that $\xl(v)=y(v)$, $\di v \in V \setminus \overline{V_0}$.
       So $y(v_0)$ is not a self-intersection point of $\xl$,
       and $f(\xl) \not=y(v_0)$ for every $\lambda \in \left(0, \delta_2\right)$. We have a contradiction. \end{proof}

\noindent {\it Proof of Corollary \ref{cy}. } (a)
It is an immediate consequence of Theorem \ref{main}.

(b) The subset of $\di \kavm$ resp. $\civm$ consisting of totally real immersions is open.
When $\dimv<\dimm$, the subset of embeddings is dense in $\di C^{k+\mathrm{sgn}(\alpha)}(V, M)$ resp. $\civm$. So there is a totally real $C^{k+\mathrm{sgn}(\alpha)}/C^{\infty}$
embedding in $\Omega$.
When $\dimv=\dimm$, the conclusion of the corollary follows from Proposition \ref{normal}
and the fact that the subset of immersions with normal crossings is dense in $\di C^{k+\mathrm{sgn}(\alpha)}(V, M)$ resp. $\civm$
(see the proof of \cite[Proposition 3.2]{gg}, where we replace the multijet transversality theorem
by  \cite[Theorem 11.2.2]{ma}, so the arguments still hold if $V$
is a $C^{\infty}$ manifold with or without boundary).

(c) If $\di j=2, 3, \cdots$ and
    $\di \dim_{\bc} M \ge \left[\frac{3}{2} \dim_{\br} V \right]$,
    then the subset of totally real immersions is dense in $C^{j}(V, M)$ (for the idea of a proof
    see \cite[Appendix]{ja}, where we replace
    the jet transversality theorem by \cite[Theorem 11.1.5]{ma}, so similar arguments
    can be applied to the general case discussed here). The conclusion of the corollary follows from (b).
    \qed

    A range decreasing holomorphic map is not always a pullback operator. Let $\Omega_j$ be the
    component of $\cistpo$ consisting of
    elements with topological degree $j \not=0$. Note that $x(S^2)=\bp^1$
    for all $x \in \Omega_j$. Thus every holomorphic map $\Omega_j
    \to \cistpo$ is range decreasing.
    Any $\gamma \in PGL(2, \bc) \setminus \{1\}$ can be considered as a
    holomorphic automorphism of $C^{\infty}(S^2, \py)$. It is straightforward
    to verify that $\gamma|_{\Omega_j}$ is not (the restriction of) a pullback operator.

%%%%%%%%%%%%%%%%%%%%%%%%%%%%%%%%%%%%%%%%%%%%%%%%%%%%%%%%%%%%%%%%%%%%%%%%%%%%%%%%%%%%%%%%%%%%%%%%%%%%%%%%%%%%%%%%%%%%%%%%%%%%%%%%
\section{Range decreasing group homomorphisms \label{group}}

In this section, we give (partial) answers to Question 1. In particular, we prove Theorem \ref{linear}
and Corollary \ref{sz}. Consider $\fvrn$ as an additive group. Let $G$ be an additive group and
$\di \varphi: \fvrn \to G$ a group homomorphism. We write
$\cs_{\varphi}$ for the subset of $V$ consisting of points $v$ with the following property:
for any open neighborhood $\fv$ of $v$, there exists $x \in \fvrn$ such that $\supp x \subset \fv$ and $\varphi(x) \not=0 \in G$.

\begin{prop} \label{sg}
  If $\varphi \not\equiv 0$, then $\di \cs_{\varphi} \not= \emptyset$.
\end{prop}
\begin{proof}
  Suppose that $\cs_{\varphi}=\emptyset$. Then for any $v \in V$, there exists an open neighborhood $\fv_v$ of $v$ such that
  $\varphi(x)=0$ for any $x \in \fvrn$ with $\supp x \subset \fv_v$. Choose a finite subcover
  $\{\fv_{v_j}\}$ of the open cover $\di \{\fv_v: v \in V\}$ of $V$,
  and take a $C^{\infty}$ partition of unity $\di \{\chi_j \}$ subordinate to $\{\fv_{v_j}\}$.
  Then $ \varphi(x)=\sum_j \varphi\left(\chi_j x \right)=0$ for all $x \in \fvrn$, this gives a contradiction.
\end{proof}

\begin{prop} \label{rn}
 Let $f: \fovrn \to \ftwrn$, where $n\ge 2$, be a real linear range decreasing map.
  Then there exists $\phi \in \tovtw$ such that $\di f=\phi^{\ast}$.
\end{prop}
\begin{proof}
For any $w \in W$, $E_w \circ f: \fovrn \to \br^n$ is a real linear range decreasing map.
It is enough to assume that $\dim W=0$ (see (\ref{evw})).

Let $e_1, \cdots, e_n$ be the standard basis
of $\br^n \subset \fovrn$. For any $x \in \fovrn$, we write
$x=\sum_{j=1}^n x_j e_j$, where $x_j \in \fovr$. In view of
the range decreasing condition, we have
$$f(x)=\sum_{j=1}^n f_j(x_j) e_j, $$
where $f_j: \fovr \to \br$ is real linear and range decreasing, $j=1, \cdots, n$.
Note that $f_j|_{\br}=\mathrm{id}$. Take $\vb \in \cs_{f_1}$ (see Proposition \ref{sg}). Next we show that
\begin{equation} \label{vw}
  \cs_{f_1}=\cdots=\cs_{f_n}=\{\vb\}.
\end{equation}
Let $v_j \in \cs_{f_j}$, where $j \not=1$. Assume that $v_j \not=\vb$. Then we can find
an open neighborhood $\mathfrak{V}$ of $\vb$,
an open neighborhood $\mathfrak{V}_j$ of $v_j$ with $\mathfrak{V} \cap \mathfrak{V}_j=\emptyset$ and
 $x_1, x_j \in \fovr$ such that $\supp x_1 \subset \mathfrak{V}$, $\supp x_j \subset \mathfrak{V}_j$,
$f_1(x_1) \not=0$ and $f_j(x_j) \not=0$. By the range decreasing condition, for any $t \in (0, 1)$,
there exists $v_t \in V$ such that
\begin{eqnarray}
  f\left(tx_1 e_1+(1-t)x_j e_j \right) &=& tf_1(x_1)e_1+(1-t) f_j(x_j) e_j \hspace{2mm} \label{e1ej}
  \\ &=& tx_1(v_t) e_1+(1-t) x_j(v_t) e_j. \label{vt}
\end{eqnarray}
Since $\supp x_1 \cap \supp x_j=\emptyset$, the value of (\ref{vt}) is either in $\{\br e_1\}$ or
in $\{\br e_j \}$. On the other hand, the value of (\ref{e1ej}) is not contained in $\{\br e_1\} \cup \{\br e_j \}$. We have
a contradiction. Therefore $v_j=\vb$ and
$\cs_{f_j}=\{\vb\}$. Similarly $\cs_{f_1}=\{\vb\}$.

For any $x_1 \in \fovr$ with $x_1(\vb) \not=0$, we claim that
\begin{equation} \label{notzero}
f_j(x_1) \not=0, \hspace{2mm} j=1, \cdots, n.
\end{equation}
Otherwise $f_{j_1}(x_1)=0$ for some $1 \le j_1 \le n$.
Take an open neighborhood $\mathfrak{V}_1$ of $\vb$ such that $x_1(v) \not=0$ for any $v \in \mathfrak{V}_1$.
If $j_2 \not=j_1$, by (\ref{vw}), we can find $x_2 \in \fovr$ with $\supp x_2 \subset \mathfrak{V}_1$ and
$f_{j_2}(x_2) \not=0$. Note that
\begin{equation*} \label{ewf}
f(x_1 e_{j_1}+x_2 e_{j_2})=f_{j_2}(x_2) e_{j_2} \in \{\br e_{j_2}\} \setminus \{\zero\}.
\end{equation*}
As $f$ is range decreasing, there exists $v_0 \in V$ such that $f(x_1 e_{j_1}+x_2 e_{j_2})$ is equal to
\begin{equation*}
  (x_1 e_{j_1}+x_2 e_{j_2})(v_0)=\left\{\begin{array}{ll}
    x_1(v_0) e_{j_1}, & v_0 \not\in \mathfrak{V}_1, \\
    x_1(v_0) e_{j_1}+x_2(v_0) e_{j_2}, & v_0 \in \mathfrak{V}_1,
  \end{array} \right.
\end{equation*}
which is not in $\{\br e_{j_2}\} \setminus \{\zero\}$. This gives a contradiction.

Note that $f|_{\br^n}=E_{\vb}|_{\br^n}=\mathrm{id}$.
So both $\ker f$ and $\ker E_{\vb}$
are of codimension $n$. It follows from (\ref{notzero}) that $\ker f \subset \ker E_{\vb}$.  Thus $\ker f=\ker E_{\vb}$
and $f=E_{\vb}$.
\end{proof}

\begin{prop} \label{lgh}
  Any range decreasing group homomorphism $f:$ $\fovgz$ $\to$ $\ftwg$ resp. $f: \fovg \to \ftwg$, where $\dim_{\br} \cg=1, 2, \cdots$,
  is a Lie group homomorphism.
\end{prop}
\begin{proof} The inclusion $\ci : \ftwg \to \cwg$ is an injective Lie group homomorphism.
It follows from the range decreasing condition that the group homomorphism $\ci \circ f$ is continuous at $1 \in \fovg$.
So it is $C^{\infty}$ (see \cite[Theorem IV.1.18]{towards}).
Let $\exp: \fg \to \cg$ be the exponential map of $\cg$, $\cb$ a ball centered at $\zero \in \fg$
such that $\exp|_{\cb}: \cb \to \exp(\cb)$ is a diffeomorphism
and $\cih$ the Lie algebra homomorphism induced by $\ci$. Recall
that the exponential map of $\fvg$ is simply $\di \exp_{\ast}^{\cf, V}$ (see Subsection \ref{family}).
Define
$$\di \fh(\xh)=\left( \exp_{\ast}^{\cf_2, W} \right)^{-1} \circ f \circ \exp_{\ast}^{\cf_1, V}(\xh) \in \ft(W, \cb), \hspace{2mm} \xh \in \fo(V, \cb).$$
 Then
\begin{equation*}
  \cih \circ \fh(\xh)=\left(\exp_{\ast}^{C, W} \right)^{-1} \circ \ci \circ f \circ \exp_{\ast}^{\cf_1, V}(\xh) \in C(W, \cb), \hspace{2mm} \xh \in \fo(V, \cb).
\end{equation*}
Thus $\cih \circ \fh$ is (the restriction of) the Lie algebra homomorphism induced by $\ci \circ f$. Now $
\di \fh=(\cih)^{-1} \circ (\cih \circ \fh)$ is real linear. Hence it is well defined
on all of $\fovfg$. Since $f$ is range decreasing, we have
\begin{equation} \label{fhr}
  \fh(\xh)(W) \subset \xh(V), \hspace{2mm} \xh \in \fovfg.
\end{equation}
It follows from Closed Graph Theorem that $\fh$ is continuous. So $f$ is $C^{\infty}$.
\end{proof}

\noindent {\it Proof of Theorem \ref{linear}.} (i) By (\ref{fhr}), Propositions \ref{rn} and \ref{lgh},
$f$ is a Lie group homomorphism which induces the same Lie algebra homomorphism as a pullback operator $\pa=\pa_{\cg}$.
So $f=\pa$ on $\fovgz$.

(ii) By (i), the closed subset $\di \{x \in \fovg: f(x)=\pa(x) \}$ is also open. So $f(x_0)=\pa(x_0)$ implies that
$f=\pa$ on the component of $\fovg$ containing $x_0$.

Let $x_1 \in \fovg$ be such that $x_1(V)$ is nowhere dense in $\cg$. We claim that $f(x_1)=\pa(x_1)$.
Otherwise there exists $w_1 \in W$ with $f(x_1)(w_1) \not=
x_1\left(\phi(w_1) \right)$.
Note that $\di f(x_1)(w_1)$, $x_1\left(\phi(w_1) \right)$ $\in$ $x_1(V)$,
hence they are in the same component of $\cg$. Fix a left invariant metric on $\cg$ and set
$$\di \delta=\dist\left(f(x_1)(w_1), x_1\left(\phi(w_1) \right) \right)>0.$$
Choose an open neighborhood $\tilde{U}$ of $ x_1\left(\phi(w_1) \right) \in \cg$ with diameter $\di <\delta/3$.
Then there exist an open contractible neighborhood $U$ of $1 \in \cg$ and an open neighborhood $V_0$ of $\phi(w_1) \in V$ such that
\begin{eqnarray*} \label{2d}
  & \dist \left(f(x_1)(w_1) a, x_1\left(\phi(w_1) \right) \right)>2\delta/3, \hspace{2mm} a \in U, \hspace{2mm} \text{and} & \\
 & \label{1d} x_1(v) a \in \tilde{U}, \hspace{2mm} v \in V_0, \,\,a \in U. &
\end{eqnarray*}
The triangle inequality gives us
\begin{equation} \label{gdelta}
  \dist \left(f(x_1)(w_1)x(\phi(w_1)), x_1(v)x(v)  \right)>\delta/3, \hspace{2mm} v \in V_0, \hspace{2mm} x \in \fo\left(V, U\right).
\end{equation}
Now take $x \in \fo\left(V, U\right)$ such that $x(v)=1$ for any $v \in V \setminus V_0$.
It follows from (i), (\ref{gdelta}) and the range decreasing condition of $f$ that
\begin{equation} \label{vvz}
f(x_1)(w_1)x(\phi(w_1))=f(x_1x)(w_1) \in x_1x(V \setminus V_0)=x_1(V \setminus V_0).
\end{equation}
With different choices of $x$, the left hand side of (\ref{vvz}) cannot always be contained
in the nowhere dense subset $x_1(V \setminus V_0)$, we have a
contradiction. So $f(x_1)=\pa(x_1)$.

If $\dimv<\dim_{\br} \cg$,
then the image of any $C^1$ map $V \to \cg$ is nowhere dense. As any component of $\fovg$
contains a $C^1$ element, we have $f=\pa$ on $\fovg$. \qed

Theorem \ref{linear} and Proposition \ref{rn} do not hold if $\dim_{\br} \cg=1$.
Let $v_1, v_2 \in V$ be two different points and $g: W \to (0, 1)$ an $\ft$ map. Define a real linear map
$f_0: \fovr \to \ftwr$ by
$$f_0(x)(w)=g(w)x(v_1)+\left(1-g(w) \right) x(v_2), \hspace{1mm} w \in W, \hspace{1mm} x \in \fovr.$$
It is clear that $f_0$ is range decreasing. We claim that $f_0$ is not a pullback operator.
Assume that $f_0=\pa$ for some $\phi: W \to V$. Fix a point $w_0 \in W$. Choose $x_1 \in \fovr$
with $x_1(v_1)=0$ and $x_1(v_2) \not=0$. Then $f_0(x_1)(w_0) \not=x_1(v_1)$, hence $\phi(w_0) \not=v_1$. Similarly $\phi(w_0) \not=v_2$.
Take $x_2 \in \fovr$ with $x_2(v_1)=x_2(v_2)=0$ and $x_2\left(\phi(w_0)\right) \not=0$. Then $f_0(x_2)(w_0) \not=\pa(x_2)(w_0)$,
which is a contradiction.

\noindent {\it Proof of Corollary \ref{sz}.} One direction being trivial, we shall only verify the  necessity part
of the claim.

(i) First we consider the special case when $\dim W=0$ (i.e. $\ftwg=\cg$).
It follows from (\ref{nonzero}) that $\di f|_{\cg}: \cg \to \cg$ is injective, which implies that
$f|_{\cg} \in Aut(\cg)$. Let $h=\left(f|_{\cg} \right)^{-1} \circ f$. Then $h|_{\cg}=\mathrm{id}$
and $h$ still satisfies (\ref{nonzero}). For any $x \in \fovgz$, we have $\di h\left(x (h(x))^{-1} \right)=1$,
where $\di (h(x))^{-1} \in \cg \subset \fovgz$.
By (\ref{nonzero}), $1 \in x (h(x))^{-1}(V)$.  Hence $h(x) \in x(V)$ (i.e. $h$ is range decreasing).
By Theorem \ref{linear}, $\di h=E_v$ for some $v \in V$. So $\di f=f|_{\cg} \circ E_v$.

(ii) For the general case, it follows from (i) that there exists a map $\phi: W \to V$ such that
$$E_w \circ f=(E_w \circ f)|_{\cg} \circ E_{\phi(w)}.$$
Define
$$\di \gf: W \ni w \mapsto (E_w \circ f)|_{\cg} \in Aut(\cg).$$
For any $a \in \cg$, the map $\di W \ni w \mapsto \gf(w)(a) \in \cg$ is just $f(a) \in \ftwg$.
Thus $\di W \ni w \mapsto d_1 \gf(w) \in Aut(\fg)$  is an $\ft$ map,
where $Aut(\fg)$ is the automorphism group of the Lie algebra $\fg$ of $\cg$.
Recall that the map $Aut(\cg) \ni \gamma \mapsto d_1 \gamma \in Aut(\fg)$
is an injective Lie group homomorphism onto a closed subgroup of $Aut(\fg)$ (e.g. see \cite[Subsection 11.3.1]{hn}).
So $\di \gf \in \ftwag$. Now $\pa=\left(\gf \right)^{-1} \circ f$ is a well defined map from $\di \fovgz$
to $\ftwg$. By Proposition \ref{fp}, we have $\phi \in \tovtw$. \qed

%%%%%%%%%%%%%%%%%%%%%%%%%%%%%%%%%%%%%%%%%%%%%%%%%%%%%%%%%%%%%%%%%%%%%%%%%%%%%%%%%%%%%%%%%%%%%%%%%%%%%%%%%%%%%%%%%%%%%%%%%%%%%%%%

\section{Decomposition of holomorphic maps \label{nm}}

In this section, we answer Question 3 partially. In particular, we prove Theorem \ref{kernel} and Corollary \ref{iff}.

\begin{lemma} \label{tg}
  Let $\sm$ be a connected complex locally convex manifold,
  $\cg$ a positive dimensional connected complex Lie group and $\di g: \sm \times \fovgz \to \ftwg$ a holomorphic map
  such that for any $\hax \in \sm$, $\di g(\hax, \cdot)$ is a group homomorphism with
  $\di g\left(\hax, \fovgo \right) \subset \ftwgo.$
  Then there exist $\phi \in \tovtw$ independent of $\hax$ and $\gamma_{g(\hax, \cdot)}$ $\in$ $\ftwag$ with
  \begin{equation*}
    g(\hax, \cdot)=\gamma_{g(\hax, \cdot)} \circ \pa_{\cg}, \hspace{2mm} \hax \in \sm
  \end{equation*}
  (cf. (\ref{ghd})).
\end{lemma}
\begin{proof}
  The group $Aut(\cg)$ of complex Lie group automorphisms of $\cg$ is a finite dimensional complex Lie group (e.g. see \cite[Section 15.4]{hn}).
  For any fixed $w \in W$, the map $\di \sm \ni \hax \mapsto E_w \circ g(\hax, \cdot)|_{\cg} \in Aut(\cg)$
  is holomorphic. It follows from Corollary \ref{sz} that for any $\hax \in \sm$, there exists $\phi_{\hax} \in \tovtw$ such that
  $$\di \left(E_w \circ g(\hax, \cdot)|_{\cg} \right)^{-1} \circ \left(E_w \circ g(\hax, \cdot) \right)=E_{\phi_{\hax}(w)}: \fovgz \to \cg.$$
  Therefore for any $\tx \in \fovgz$, the map $\sm \ni \hax \mapsto E_{\phi_{\hax}(w)}(\tx) \in \cg$ is holomorphic.
  Let $\rho: \br \to \cg$ be a $C^{\infty}$ embedding and $\mu \in \fovr$. Then $\rho \circ \mu \in \fovgz$,
  and $\sm \ni \hax \mapsto E_{\phi_{\hax}(w)}\left(\rho \circ \mu \right) \in \cg$ is a holomorphic map whose image is contained in $\rho(\br)$.
  Thus $\di E_{\phi_{\hax}(w)}\left(\rho \circ \mu \right)$ is independent of $\hax$ for every $\mu \in \fovr$,
  which implies that $\phi_{\hax}(w)$ is independent of $\hax$.
\end{proof}

\noindent {\it Proof of Theorem \ref{kernel}.} It is enough to show that for any $y \in \Omega$, there exist a connected open neighborhood
$\su_y \subset \Omega$ of $y$ and $v=v(y) \in V$ independent of $x \in \su_y$ such that $\ker d_x f=\ker E_{v, \fovcn}$ for every $x \in \su_y$.
As $\Omega$ is connected, $v$ is also independent of $y$.

Let $\Phi: U \to \bc^n$ be a local chart of $M$ with $f(y) \in U$ and $\Phi\left(f(y) \right)=\zero$,
$\varphi_y: \su_y \to \fovcn$ a local chart of $\fovn$ as in
\cite[Subsection 1.1]{ls}
such that $\su_y$ is connected, $y \in \su_y$, $\varphi_y(y)=\zero$ and $\di f\left(\su_y \right) \subset U$,
$$\di h=\Phi \circ f \circ \varphi_y^{-1}: \varphi_y(\su_y) \to \bc^n$$
and $\hax_1, \hax_2 \in \varphi_y(\su_y)$ with $\hax_1(v) \not= \hax_2(v)$, $v \in V$.  By \cite[Lemma 1.1]{ls},
\begin{equation} \label{distinct}
\varphi_y^{-1}(\hax_1)(v) \not= \varphi_y^{-1}(\hax_2)(v), \hspace{1mm}\hspace{1mm} v \in V.
\end{equation}
It follows from (\ref{x1x2}) that $h(\hax_1) \not=h(\hax_2)$. In particular, $h|_{\varphi_y(\su_y) \cap \bc^n}$ is injective, so it is an embedding.

Next we show that $\di d_{\zero} h \left(\fovcno \right) \subset \bc^n \setminus \{\zero\}$.
Let $\di \lambda \in \bc \setminus \{0\}$, $\tx \in \fovcno$, $\di A=d_{\zero}\left(h|_{\varphi_y(\su_y) \cap \bc^n}\right) \in GL(\bc^n)$,
$$\di m(\tx)=\min_{v \in V} \left|\tx(v) \right|>0, \hspace{2mm} |A|=\min_{\zeta \in \bc^n, |\zeta|=1} |A (\zeta)|>0$$
and $B \subset \varphi_y(\su_y) \cap \bc^n$ a small ball centered at $\zero$ such that
$$\left|h(\zeta) \right|>|A||\zeta|/2, \hspace{2mm} \zeta \in B \setminus \{\zero\}.$$
 Note that $h(\lambda \tx) \in h(B)$ when $|\lambda|$ is small enough.
In this case, there is a unique $\zeta(\lambda) \in B$ with $h(\lambda \tx)=h(\zeta(\lambda))$.
Since $h$ has the property as in (\ref{x1x2}), we must have $\zeta(\lambda) \in \lambda \tx(V)$. Therefore
$$\left| h(\lambda \tx) \right| > |A| |\zeta(\lambda)|/2 \ge  |A| |\lambda| m(\tx)/2,$$
 which implies that $d_{\zero} h(\tx) \in \bc^n \setminus \{\zero\}$.
 Similarly $\di d_{\hax} h \left(\fovcno \right) \subset \bc^n \setminus \{\zero\}$ for every $\hax \in \varphi_y(\su_y)$.
Applying Lemma \ref{tg} (where we set $\dim W=0$) to the holomorphic map
$$\di \varphi_y(\su_y) \times \fovcn \ni (\hax, \tx) \mapsto d_{\hax} h (\tx) \in \bc^n,$$
we obtain that
$\ker d_{\hax} h=\ker E_{v, \fovcn}=\ker d_x f$, where $x=\varphi_y^{-1}(\hax) \in \su_y$ and $v \in V$ is independent of $\hax \in \varphi_y(\su_y)$.
\qed

Let $x_1 \in \fvpy$ be such that the bundle $x_1^{\ast} T\bp^1$ is trivial.
By (\ref{distinct}),
we can find $x_2, x_3$ in a neighborhood of $x_1$ with $x_i(v) \not= x_j(v)$ for all $v \in V$, $i, j=1, 2, 3$, $i \not=j$.
It follows from the proof of Proposition \ref{2points}(b) that
there exists an element of $\fvpglt$ which maps $0 \in \py$ to $x_1$. In particular,
the action of $\fvpglt$ on $\fvpy$
is transitive if $H^2(V, \bz)=0$. Recall that $\cf(S^3, \bp^1)$ has infinitely many components.
Therefore $\cf^0(S^3, \bp^1)$ is not an invariant subset of $\di \cf(S^3, \bp^1)$ under the action of
$\cf\left(S^3, PGL(2, \bc) \right)$.

\noindent {\it Proof of Corollary \ref{iff}.}
  Suppose that $\di f=f|_{N} \circ E_{v_0}|_{O}$ with $N_{n, f|_N} \not=\emptyset$. Take an open subset $U \subset N$
  such that $\di f|_U$ is an embedding and let $\Omega=E_{v_0}^{-1}(U) \cap O$. Then $f(\Omega)= f(U)$
  and (\ref{x1x2}) holds on $\Omega$. For the other direction, it follows from Theorem \ref{kernel}
  that $\di f=f|_{N} \circ E_{v_0}$ on an open subset of $\Omega$, thus also on the
  connected manifold $O$. \qed

Some special cases of the sufficiency part of Corollary \ref{iff} could also be proved by Theorem \ref{main} or by Corollary
\ref{cy} (instead of Theorem \ref{kernel}): By (\ref{x1x2}), $\di f|_{\Omega \cap N}$
is injective. Take
a component $U_1$ of $\Omega \cap N$. Then $\di f|_{U_1}: U_1 \to f(U_1)$ is biholomorphic.
Let $\Omega_1$ be the component of
$\Omega \cap f^{-1}\left(f(U_1) \right)$ containing $U_1$. Then $\di f(\Omega_1)=f(U_1)$.
Define $$\di \tf=\left(f|_{U_1} \right)^{-1} \circ f|_{\Omega_1}: \Omega_1 \to U_1 \subset N.$$
We only need to show that $\tf=E_{v_0}|_{\Omega_1}$ for some $v_0 \in V$.
It follows from (\ref{x1x2}) that for any $x \in \Omega_1$ and any $\zeta \in U_1 \setminus x(V)$ (if $U_1 \setminus x(V) \not=\emptyset$), we have
$\tilde{f}(x) \not=\tilde{f}(\zeta)=\zeta$.
So
$$\di \tilde{f}(x) \in x(V), \hspace{2mm} x \in \Omega_1$$
(i.e. $\tilde{f}$ is range decreasing with $\dim W=0$).
If $\tf$ satisfies the conditions of  Theorem \ref{main} or of Corollary
\ref{cy}, then we obtain $\tilde{f}=E_{v_0}|_{\Omega_1}$.

%%%%%%%%%%%%%%%%%%%%%%%%%%%%%%%%%%%%%%%%%%%%%%%%%%%%%%%%%%%%%%%%%%%%%%%%%%%%%%%%%%%%%%%%%%%%%%%%%%%%%%%%%%%%%%%%%%%%%%%%%%%%%%%%%%%%%%%%%%%%%%
\section{Holomorphic maps $\fovpnz \to \ftwpm$ \label{pn}}

In this final section, we study the decomposition of holomorphic maps $\fovpnz$ $\to$ $\ftwpm$
as in Question 3.
Let $O \subset \fovpnz$ be a connected open neighborhood of $\bp^n$ and $f: O \to \ftwpm$ a holomorphic map
with $\di \deg f=d \ge 1$ (which implies that $m \ge n$).
Then $f$ induces a map $$\di \gf: W \ni w \mapsto E_w \circ f|_{\bp^n} \in \hd.$$
For any $\zeta \in \pn$, the map $W \ni w \mapsto \gf(w)(\zeta) \in \bp^m$ is exactly $f(\zeta) \in \ftwpm$.
So
\begin{equation} \label{gf}
\gf \in \ftwhd,
\end{equation}
and $\gf$ can be considered as a holomorphic map $\ftwpn \to \ftwpm$
with $\gf|_{\bp^n}=f|_{\bp^n}$  (see Subsection \ref{family}).

\begin{lemma} \label{phi}
  Let $f, \gf$ be as above and let $\phi: W \to V$ be a map. If
    \begin{equation} \label{gfw}
     E_{w, \ftwpm} \circ f=\gf(w) \circ E_{\phi(w), \fovpnz}|_{O}, \hspace{2mm} w \in W,
     \end{equation}
     then $\phi \in \tovtw$.
\end{lemma}
\begin{proof}
It is enough to show that for any $w_0 \in W$, there is an open neighborhood $W_0$ of $w_0$
such that its closure $\overline{W_0}$ is a contractible $C^{\infty}$ submanifold of $W$, possibly with boundary, and
$\phi|_{\overline{W_0}} \in \tovtwz$.
Let $\ci: \overline{W_0} \to W$
be the inclusion and $\ci^{\ast}: \ftwpm \to \ftwzpm$ the pullback operator induced by $\ci$.
Replacing $f$ by $\ci^{\ast} \circ f$, we may assume that $W=\overline{W_0}$.
Then $\ftwpm$ is connected.
We identify the tangent space of $\ftwpm$
(resp. $\fovpnz$) at any point with $\di \ftwcm$ (resp. $\di \fovcn$).

First we consider the case when $m=n$. Shrinking $W=\overline{W_0}$ if necessary, we can find $\zeta_0 \in \bp^n \subset \ftwpn$
such that the tangent map
$\di d_{\zeta_0}\gf(w)$ is of rank $n$ for any $w \in W$.
Then the map $\kappa_{\gf, \zeta_0}: \zeta_0^{\ast}T\bp^n  \to \gf(\zeta_0)^{\ast}T\bp^n$ as in (\ref{kappa}) is an $\ft$ isomorphism
of trivial vector bundles. By (\ref{dx}), the tangent map
$\di d_{\zeta_0} \gf$ is an automorphism of the complex Banach/Fr\'echet space $\ftwcn$.
Note that $\gf(\zeta_0)=f(\zeta_0)$.
It follows from (\ref{gfw}) that
$$E_{\phi(w), \fovcn}=\left(d_{\zeta_0}\gf(w) \right)^{-1} \circ E_{w, \ftwcn} \circ d_{\zeta_0} f, \hspace{2mm} w \in W.$$
So the pullback by $\phi$ from $\fovcn$ to the space of maps $W \to \bc^n$ is exactly
$\di \left( d_{\zeta_0} \gf\right)^{-1} \circ d_{\zeta_0} f: \fovcn \to \ftwcn$.
Thus $\phi \in \tovtw$.

When $m>n$, take a hyperplane $\bp^{m-1} \subset \bp^m$ and a point $p \in \bp^m \setminus \bp^{m-1}$
such that $\gf(w)(\bp^n) \subset \bp^m \setminus \{p\}$ for any $w \in W$ (which can be done for sufficiently
small $W=\overline{W_0}$). In view of (\ref{gfw}), we have $f(O) \subset \cf_2(W, \bp^m \setminus \{p\}).$
Let $\pi: \bp^m \setminus \{p\} \to \bp^{m-1}$ be the projection
from $p$ to $\bp^{m-1}$,
$$\pi_{\ast}: \cf_2(W, \bp^m \setminus \{p\}) \ni y \mapsto \pi \circ y \in \ft(W, \bp^{m-1})$$
the induced holomorphic map and $f_1=\pi_{\ast} \circ f$. By (\ref{gfw}), we have
     \begin{equation*}
     E_{w, \ft(W, \bp^{m-1})} \circ f_1=\gamma_{f_1}(w) \circ E_{\phi(w), \fovpnz}|_{O}, \hspace{2mm} w \in W;
     \end{equation*}
and the conclusion of the lemma follows from an induction on $m$.
\end{proof}

Consider $\di \bp^n$ as $\bc^n \cup \bp^{n-1}$, where we identify $\bc^n$ with $\{ [Z_0, Z_1, \cdots, Z_n] \in \bp^n:
Z_0 \not=0\}$. Let $x \in \fvcn$ and
$$y=(y_1, \cdots, y_n) \in \cf(V, \bc^n \setminus \{\zero\}).$$
Recall the class of rational curves $\fc$ in Subsection \ref{geo}. Define
\begin{equation*} \label{cxy}
  \cxy=\{x+\lambda y \in \fvpnz: \lambda \in \bp^1 \} \in \fc,
\end{equation*}
 where $\lambda=\infty \in \bp^1$ corresponds to the point
\begin{equation} \label{ly}
\fin(y)=[0, y_1, \cdots, y_n] \in \cf(V, \bp^{n-1}) \cap \fvpnz.
\end{equation}
If $n=1$, then $\fin(y)=\infty \in \bp^1 \subset \fvpyz$.

\begin{prop} \label{not=}
  Let $f$ be as in Theorem \ref{degreeone} and assume that $\dim W=0$ (i.e. $\ftwpm=\bp^m$). Then there exist $\zeta \in \bc^n$
  and an open neighborhood $\Omega \subset \fovcn$ of $\zeta$ such that $\di f(\Omega) \subset f(\bp^n)$ and
  (\ref{x1x2}) holds on $\Omega$.
\end{prop}
\begin{proof}
  If $d=1$, then $f(\bp^n)$ is an $n$-plane in $\bp^m$.
  For any curve $C \in \fc$ in $\fovpnz$, $f|_C$ is injective and $f(C)$ is a projective line in $\bp^m$.
  Next we show that
  \begin{equation} \label{contain}
  \di f(\fovcn) \subset f(\bp^n)
  \end{equation}
  by induction on $n$. Suppose $n=1$.
 For any $x \in \fovc$, take $\zeta_1, \zeta_2 \in \bp^1 \setminus x(V)$ with $\zeta_1 \not= \zeta_2$. By Proposition \ref{2points}(b),
there is a curve $C_1 \in \fc$ through $\zeta_1, \zeta_2$ and $x$. Then $f(C_1)$ is the projective line
through $f(\zeta_1)$ and $f(\zeta_2)$. So $f(x) \in f(\bp^1)$.
Assume that (\ref{contain}) holds for $n=k-1$, where $k \ge 2$, $m \ge n$.
When $n=k$, consider $\bp^k$ as $\bc^k \cup \bc^{k-1} \cup \bp^{k-2}$.
Let $\fO \subset \cf_1(V, \bc^k \setminus \{\zero\})$ be the open subset consisting of $y$
with $\fin(y) \in \cf_1(V, \bc^{k-1})$ (thus $\di f(\fin(y)) \in f(\bp^{k-1})$).
For any $y \in \fO$, $f(C_{0, \, y})$ is the projective line through $f(0)$ and $f(\fin(y))$.
Thus $\di f\left(\fO \right) \subset f(\bp^k)$. Note that $\di \fO \cap (\bc^k \setminus \{\zero\})$
is dense in $\bc^k \setminus \{\zero\}$.
For any $x \in \fovck$,
take $\zeta \in \fO \cap (\bc^k \setminus \{\zero\})\setminus x(V)$.
By Proposition \ref{2points}(a), there is a curve $C_2 \in \fc$ through $x$ and $\zeta$.
The openness of $\fO$ implies that $C_2$ contains infinitely many points of $\fO$.
Hence $\di f(C_2) \subset f(\bp^k)$ and  (\ref{contain}) holds for $n=k$.
Let $\di \Omega=\fovcn$ and $x_1, x_2 \in \Omega$ with $x_1(v) \not= x_2(v)$ for all $v \in V$. Then
there is a curve $C_3 \in \fc$ through $x_1$ and $x_2$.
So $f(x_1) \not=f(x_2)$.

For the map $f$ in Theorem \ref{degreeone}(b),
let $\zeta_1, \zeta_2 \in \bc$ be regular points of $f|_{\bp^1}$ with $\zeta_1 \not= \zeta_2$
and $f(\zeta_1)=f(\zeta_2)$, $D \subset \bp^1$ an open neighborhood of $f(\zeta_1)$ and $U_j \subset \bc$ an open neighborhood of $\zeta_j$, $j=1, 2$,
such that $U_1 \cap U_2=\emptyset$ and $f$ maps $U_j$ biholomorphically onto $D$.
Take a connected open neighborhood $\Omega \subset \cf_1(V, U_1)$ of $\zeta_1$ with $f(\Omega) \subset D$.
Let $x_1, x_2 \in \Omega$ with $x_1(v) \not=x_2(v)$ for every $v \in V$.
Then we can find $\zeta_3=\zeta_3(x_1) \in U_2$ with $f(\zeta_3)=f(x_1)$.
By Proposition \ref{2points}(b), there is a curve $C_4 \in \fc$ through $\zeta_3$, $x_1$ and $x_2$.
As $f|_{C_4}: C_4 \to \py$ is of topological degree two, we have $f(x_1) \not= f(x_2)$.
\end{proof}

\noindent {\it Proof of Theorem \ref{degreeone}.} By Proposition \ref{not=} and Corollary \ref{iff}, there is a map $\phi: W \to V$ such that
$$E_{w} \circ f=\left(E_w \circ f \right)|_{\bp^n} \circ E_{\phi(w)}:
\fovpnz \to \bp^m, \hspace{2mm} w \in W.$$
The conclusion of the theorem follows from (\ref{gf}) and Lemma \ref{phi}. \qed

Corollary \ref{auto} follows from Theorem \ref{degreeone} and (\ref{comm}).

Next we construct holomorphic maps $\cisopo \to \ftwpm$ which are not of the form (\ref{decomposition}).
The constructions are closely related to the $C^{\infty}$ maps
\begin{equation*}
 \jk: S^1 \times \cisopo \ni (t, x) \mapsto j^k_t x \in J^k(S^1, \bp^1),
\end{equation*}
the target maps $\tau_k: J^k(S^1, \bp^1) \to \bp^1$ of the jet spaces, $k=0, 1, \cdots$, and the evaluation
\begin{equation*}
\di E=\tau_k \circ \jk: S^1 \times \cisopo \ni (t, x) \mapsto x(t) \in \bp^1.
\end{equation*}
Recall that for any fixed $t_1 \in S^1$, the map
$\di j^k_{t_1}=\jk(t_1, \cdot): \cisopo \to J^k(t_1, \bp^1)$ and the restriction of $\tau_k$ to $J^k(t_1, \bp^1)$
are holomorphic (see Section \ref{holorot}).
If $h: N \to M$ is a holomorphic map between complex manifolds and $\cl$ is a holomorphic line bundle over $M$,
then we write $\di H^0(M, \cl)$ for the space of holomorphic sections of $\cl$, and denote by $\di h^{\ast}H^0(M, \cl)$
the subspace of $\di H^0(N, h^{\ast} \cl)$ consisting of pullback sections.
For any positive integer $n_1$, the pullback of the line bundle $\co(n_1) \to \bp^1$ by
\begin{equation*} \label{eto}
E_{t_1}=\tau_k \circ j^k_{t_1}: \cisopo \to \bp^1
\end{equation*}
is actually
the bundle $\lpo$, where $\di \varphi_1=E_{t_1}^{n_1}$, in \cite[Section 4]{z03}.
By \cite[Theorem 1.2]{z03},
$n_1+1 \le \dim \hocipo <\infty.$

  For any $\nu \in \hopono$,
  $E^{\ast} \nu=(\jk)^{\ast} \tau_k^{\ast} \nu$ is a $C^{\infty}$ section of the bundle $$\di E^{\ast} \co(n_1) \to S^1 \times \cisopo,$$
  where $\tau_k^{\ast} \nu$ is a $C^{\infty}$ section of the bundle $\di \tau_k^{\ast} \co(n_1) \to J^k(S^1, \bp^1)$
  which is holomorphic on each of the submanifolds $J^k(t_1, \bp^1) \subset J^k(S^1, \bp^1)$, $t_1 \in S^1$, and
  we may consider $E_{t_1}^{\ast} \nu$ as the restriction of $E^{\ast} \nu$ to $\di \{t_1\} \times \cisopo$.
  Let $\{\sigma_a \in H^0 \left(C^{\infty}\left(S^1, \bp^1 \setminus \{a\} \right), \lpo \right): a \in \bp^1\}$ be the family of non-vanishing sections in \cite[Proposition 2.7]{z03}.
  Then there are sections $s_a \in H^0 \left(\bp^1, \co(1) \right)$ such that
  \begin{equation} \label{sigmaa}
  \sigma_a=E_{t_1}^{\ast} s_a^{n_1}|_{C^{\infty}\left(S^1, \bp^1 \setminus \{a\} \right)}
  =(j^k_{t_1} )^{\ast} \tau_k^{\ast} s_a^{n_1}|_{C^{\infty}\left(S^1, \bp^1 \setminus \{a\} \right)}, \hspace{2mm} a \in \bp^1
  \end{equation}
  (see the proof of \cite[Proposition 2.7]{z03}).  By Proposition 4.2 and Theorem 4.7 of \cite{z03},
  for any $\sigma \in \hocipo$, $\sigma/\sigma_{\infty}: \cisoc \to \bc$ is a polynomial
  in $\di x(t_1)$, $x'(t_1)$, $\cdots$, $x^{(n_1-1)}(t_1)$, where $\di x \in \cisoc$.
  Hence there is a holomorphic function $\varsigma_{\infty}$
  on $\di \tau_{n_1-1}^{-1}(\bc) \cap J^{n_1-1}(t_1, \bp^1)$ with
  $\sigma/\sigma_{\infty}=(j_{t_1}^{n_1-1})^{\ast} \varsigma_{\infty}$.
  For any $\gamma \in PGL(2, \bc)$, we have $\gamma^{\ast} \lpo \simeq \lpo$ (see \cite[Section 2]{z03}).
  Similar to  $\varsigma_{\infty}$, we can find holomorphic functions
  $\varsigma_a$ on $\di \tau_{n_1-1}^{-1}(\bp^1 \setminus \{a\}) \cap J^{n_1-1}(t_1, \bp^1)$ with
  \begin{equation} \label{ssigmaa}
  \sigma/\sigma_{a}=(j_{t_1}^{n_1-1})^{\ast} \varsigma_{a}, \hspace{2mm} a \in \bp^1.
  \end{equation}
  Combination of (\ref{sigmaa}) and (\ref{ssigmaa}) gives
  \begin{equation} \label{nomo}
  \hocipo=\jahojno.
  \end{equation}
   If $n_1 \ge 2$, then
   $$\jahojkm \subsetneqq \jahojk,$$
    where $k=1, \cdots, n_1-1$, and
    $$\left(j_{t_{1}}^{0}\right)^{*} H^{0}\left(J^{0}\left(t_{1}, \mathbb{P}^{1}\right), \tau_{0}^{*} \mathcal{O}\left(n_{1}\right)\right)
    =\ahopono.$$
    If we consider $t_1$ in the expression of $\di \sigma/\sigma_{\infty}$ as a variable in $S^1$, then we obtain
  a $C^{\infty}$ function $S^1 \times \cisoc \to \bc$ which is the pullback
  of a $C^{\infty}$ function $\tilde{\varsigma}_{\infty}$ on the open subset $\di \tau_{n_1-1}^{-1}(\bc)$ of $J^{n_1-1}(S^1, \bp^1)$ by $\jnom$.
  Similarly, each $\di \sigma/\sigma_{a}$ induces a $C^{\infty}$ function $\tilde{\varsigma}_{a}$
  on the open subset $\di \tau_{n_1-1}^{-1}(\bp^1 \setminus \{a\})$ of $J^{n_1-1}(S^1, \bp^1)$.
  Hence we obtain a $C^{\infty}$ section $\tilde{\sigma}$ of the bundle $\di E^{\ast} \co(n_1) \to S^1 \times \cisopo$
  such that
  $$\left(\tilde{\sigma}/E^{\ast} s_{a}^{n_1} \right)|_{S^1 \times \cisopoa}=\left(\jnom \right)^{\ast} \tilde{\varsigma}_{a}$$
  (and $\sigma$ can be considered as the restriction of $\tilde{\sigma}$ to $\di \{t_1\} \times \cisopo$). Let
  $$\fepo=\left\{\tilde{\sigma}: \sigma \in \hocipo \right\}.$$

  Choose sections $\di \ts_1, \cdots, \ts_{m'} \in \ehopono \subset \fepo$, where $m' \ge 2$,
  without common zeros
  and take linearly independent sections
  $$\ts_{m'+1}, \cdots, \ts_{m+1} \in \fepo \setminus \ehopono$$
  (where $\fepo \setminus \ehopono \not=\emptyset$ if $n_1 \ge 2$).
  Then the subspace of $\fepo$ spanned by $\ts_1, \cdots, \ts_{m+1}$ gives rise to
  a $C^{\infty}$ map $\sg: S^1 \times \cisopo \to \bp^m$ (e.g. see \cite[Section 1.4]{gh}) such that
  $\sg(t_1, \cdot): \cisopo \to \bp^m$ is a holomorphic map of degree $n_1$ for every $t_1 \in S^1$.
  The map $\sg$ induces a holomorphic map
  \begin{equation} \label{gammao}
  \di \tilde{g}: \cisopo \to \cisopm, \hspace{2mm} \text{where} \hspace{2mm} \tilde{g}(x)(t_1)=\sg(t_1, x).
  \end{equation}
  By (\ref{nomo}), we may consider $E_{t_1} \circ \tilde{g}=\sg(t_1, \cdot)$ as the composition of $\di j^{n_1-1}_{t_1}$
  and a holomorphic map $\di J^{n_1-1}(t_1, \bp^1) \to \bp^m$.
  If $n_1 \ge 2$ and $m \ge m'$, then $\tilde{g}$ is not of the form (\ref{decomposition}).

More generally, let $\di \Upsilon: \bp^{m_1} \times \cdots \times \bp^{m_p} \to \bp^m$, $g_j: \fovpnz \to \fovpmj$
be holomorphic maps and $\phi_j \in \tovtw$, $j=1, \cdots, p$. Recall that $\Upsilon$ induces
a holomorphic map
$$\Upsilon_{\ast}: \ftwpmo \times \cdots \times \ftwpmp \to \ftwpm$$
as in (\ref{push}).  Define a holomorphic map
\begin{eqnarray}
& \di g=g_{\Upsilon, g_1, \cdots, g_p, \phi_1, \cdots, \phi_p}: \fovpnz \to \ftwpm \hspace{2mm} \text{by} & \notag \\ &
 \label{general}  g(x)=\Upsilon_{\ast} \left((\phi_1)^{\ast}_{\bp^{m_1}}\left(g_1(x)\right), \cdots, (\phi_p)^{\ast}_{\bp^{m_p}}\left(g_p(x)\right) \right). &
\end{eqnarray}
Now choose $g_j$ to be maps $\cisopo \to \cisopmj$ as in (\ref{gammao}), $\phi_j \in \tcisotw$, where $j=1, \cdots, p$,
and take $\Upsilon$ to be injective. If  $p \ge 2$ and not all $\phi_1, \cdots, \phi_p$ are the same,
then the map $g: \cisopo \to \ftwpm$ as in (\ref{general}) is not of the form (\ref{decomposition}).

By similar arguments as above, we could also construct holomorphic maps $\di f: \fosopo \to \ftwpm$,
where $\fo=C^k$
or $W^{k, p}$  (see \cite[Section 4]{z03}).


\begin{thebibliography}{bhs}
\bibitem{bhs} B. Bojarski, P. Haj{\l}asz, P. Strzelecki, Sard's theorem for mappings in H\"older and Sobolev
spaces, manuscripta math. 118 (2005), 383-397.


\bibitem{d} S. Dineen, Complex analysis on infinite dimensional
spaces, Springer, London, 1999.

\bibitem{fe} S. Feder, Immersions and embeddings in complex projective spaces, Topology 4 (1965), 143-158.

\bibitem{fo} F. Forstneri\v{c}, Stein manifolds and holomorphic mappings, Springer, Berlin, 2011.

\bibitem{gh} P. Griffiths, J. Harris, Principles of Algebraic
Geometry, John Wiley \& Sons, New York, 1978.

\bibitem{gg} M. Golubitsky, V. Guillemin, Stable mappings and their
singularities, Springer, New York, 1973.


\bibitem{rh} R. Hamilton, The inverse function theorem of Nash and Moser, Bull. Amer.
Math. Soc. 7 (1982), 65-222.

\bibitem{h} M. Herv\'{e}, Analyticity in infinite dimensional
spaces, Walter de Gruyter, Berlin, 1989.

\bibitem{hn} J. Hilgert, K.-H. Neeb, Structure and geometry of Lie groups, Springer, New York, 2012.

\bibitem{ho} G. Hochschild, The automorphism group of a Lie group, Trans. Amer. Math. Soc. 72 (1952), 209-216.

\bibitem{ja} P. T. Ho, H. Jacobowitz, P. Landweber, Optimality for totally real immersions and
independent mappings of manifolds into $\bc^N$, New York J. Math. 18 (2012), 463-477.

\bibitem{kr} N. Krikorian, Differential structures on
function spaces, Trans. Amer. Math. Soc. 171 (1972), 67-82.

\bibitem{l04} L. Lempert, Holomorphic functions on
(generalised) loop spaces, Math. Proc. R. Ir. Acad. 104A (2004),
35-46.

\bibitem{ls} L. Lempert, E. Szab\'o, Rationally connected
varieties and loop spaces, Asian J. Math. 11 (2007), 485-496.

\bibitem{lz} L. Lempert, N. Zhang, Dolbeault cohomology of a
loop space, Acta Math. 193 (2004), 241-268.

\bibitem{ma} J. Margalef-Roig, E. O. Dominguez, Differential topology, North-Holland, Amsterdam, 1992.

\bibitem{mi} P. W. Michor, Manifolds of differentiable mappings, Shiva Publishing Limited, Orpington, 1980.

\bibitem{mo} J. Mostovoy, Spaces of rational maps and the Stone-Weierstrass theorem, Topology 45 (2006), 281-293. 

\bibitem{towards} K.-H. Neeb, Towards a Lie theory of locally convex groups, Japan. J. Math. 1 (2006), 291-468. 

\bibitem{ns} K.-H. Neeb, H. Sepp\"anen, Borel-Weil theory for groups over
commutative Banach algebras, J. reine angew. Math. 655 (2011), 165-187. 

\bibitem{pa} R. S. Palais, Foundations of global
non-linear analysis, Benjamin, New York, 1968.

\bibitem{ps} A. Pressley, G. Segal, Loop groups, Oxford
University Press, New York, 1986.

\bibitem{ya} K. Yamaguchi, Fundamental groups of spaces of holomorphic
maps and group actions, J. Math. Kyoto Univ. 44 (2004), 479-492. 

\bibitem{z03} N. Zhang, Holomorphic line bundles on the loop
space of the Riemann sphere, J. Differ. Geom. 65 (2003), 1-17.

\bibitem{z10} N. Zhang, The Picard group of the loop
space of the Riemann sphere, Int. J. Math. 21 (2010), 1387-1399.

\bibitem{z17} N. Zhang, Holomorphic automorphisms of the loop space of
$\mathbb{P}^n$, Commun. Anal. Geom. 25 (2017), 709-718. 
\end{thebibliography}
\end{document}